\documentstyle[12pt]{article}
\author{S. Shelah
\thanks{The first author was supported by the United States Israel Binational
	Science Foundation, Publication 370}
	\\Institute of Mathematics\\Hebrew University, Jerusalem \and L. Soukup
	\thanks{The second author was supported by the
	Hungarian National Foundation for Scientific Research grant no. l805}
	\\Mathematical Institute of the\\Hungarian Academy of Sciences}
	\title{On the number of  non-isomorphic
	subgraphs}
	\begin{document}
	\pagestyle{myheadings}\markright{{\rm\today}}
	\maketitle

	\newtheorem{theorem}{Theorem}[section]
	\newtheorem{Definition}[theorem]{Definition}
	\newtheorem{proposition}[theorem]{Proposition}
	\newtheorem{lemma}[theorem]{Lemma}
	\newtheorem{sublemma}{Sublemma}[theorem]

	\def\ujbekezdes{\par}
	\newenvironment{definition}
	{\begin{Definition}\rm}{\end{Definition}}
	\newenvironment{demo}[1]
	{{\ujbekezdes}\noindent{\bf #1:}\ \rm}{{\ujbekezdes}}
	\newenvironment{case}[2]
	{{\ujbekezdes}\noindent{\bf Case\ #1:}{\em #2}\ujbekezdes \rm}{{\ujbekezdes}}

	\newcommand{\proof}{\begin{demo}{Proof}}
	\newcommand{\Proof}[1]{\begin{demo}{Proof of #1}}
	\newcommand{\eproof}{\ \rule{2mm}{2mm}\end{demo}\vspace{1mm}}

	\newcommand{\arablabel}{
	          \renewcommand{\labelenumi}{{\rm (\arabic{enumi})}}
	          \renewcommand{\theenumi}{{\rm (\arabic{enumi})}}
	                    }

	\newcommand{\alabel}{
	          \renewcommand{\labelenumi}{{\rm (\alph{enumi})}}
	          \renewcommand{\theenumi}{{\rm (\alph{enumi})}}
	                    }
	\newcommand{\Alabel}{
	          \renewcommand{\labelenumi}{{\rm (\Alph{enumi})}}
	          \renewcommand{\theenumi}{{\rm (\Alph{enumi})}}
	                    }
	\newcommand{\rlabel}{
	          \renewcommand{\labelenumi}{{\rm (\roman{enumi})}}
	          \renewcommand{\theenumi}{{\rm (\roman{enumi})}}
	                    }
	\newcommand{\Rlabel}{
	         \renewcommand{\labelenumi}{{\rm(\Roman{enumi})}}
	         \renewcommand{\theenumi}{{\rm(\Roman{enumi})}}
	                    }

	\newcommand{\litem}{\renewcommand{\labelitemi}{{\bf --}}}

	\newdimen\elso
	\newdimen\masodik
	\newdimen\harmadik
	\def\tagg#1#2$${
	\settowidth{\elso}{\hbox{$\displaystyle #2$}}
	\settowidth{\masodik}{\hbox{$\displaystyle (#1)$}}
	\harmadik \textwidth
	\advance\harmadik by -\elso
	\advance\harmadik by -\masodik
	\advance\harmadik by -\masodik
	\divide\harmadik by 2
	\hbox{$\displaystyle (#1)$}
	\makebox[\harmadik]{}
	\hbox{$\displaystyle #2$}
	\makebox[\harmadik]{}
	\makebox[\masodik]{}
	$$}

	\def\multline{\begin{array}{c}}
	\def\endmultline{\end{array}}

	\newcommand{\limo}{\mbox{\rm {Lim}}}
	\newcommand{\leo}{<_{\mbox{\rm{On}}}}
	\newcommand{\oo}{{{\omega}_1}}
	\newcommand{\acal}{{\cal A}}
	\newcommand{\bcal}{{\cal B}}
	\newcommand{\ccal}{{\cal C}}
	\newcommand{\dcal}{{\cal D}}
	\newcommand{\fcal}{{\cal F}}
	\newcommand{\gcal}{{\cal G}}
	\newcommand{\hcal}{{\cal H}}
	\newcommand{\ical}{{\cal I}}
	\newcommand{\jcal}{{\cal J}}
	\newcommand{\kcal}{{\cal K}}
	\newcommand{\mcal}{{\cal M}}
	\newcommand{\ncal}{{\cal N}}
	\newcommand{\pcal}{{\cal P}}
	\newcommand{\qcal}{{\cal Q}}
	\newcommand{\scal}{{\cal S}}
	\newcommand{\tcal}{{\cal T}}
	\newcommand{\wcal}{{\cal W}}
	\newcommand{\atil}{\tilde a}
	\newcommand{\btil}{\tilde b}
	\newcommand{\adot}{\dot \alpha}
	\newcommand{\bdot}{\dot \beta}
	\newcommand{\abar}{\bar a}
	\newcommand{\hhbar}{\dot h}
	\newcommand{\fdot}{\dot f}
	\newcommand{\bbar}{\bar b}
	\def\rest{{\tsize |}\raise 5pt \hbox{$@!@!\!\sssize\backslash$}}
	\newcommand{\ahat}{\hat a}
	\newcommand{\bhat}{\hat b}
	\newcommand{\chat}{\hat c}
	\newcommand{\dhat}{\hat d}
	\newcommand{\dplus}{{\diamondsuit^{\scriptscriptstyle +}}}
	\newcommand{\fini}{{<{\omega}}}
	\newcommand{\setm}{\setminus}
	\newcommand{\empt}{\emptyset}
	\newcommand{\subs}{\subset}
	\newcommand{\bijp}{\mbox{\rm {Bij}${}_p$}}
	\newcommand{\isop}{\mbox{\rm{Iso}${}_p$}}
	\newcommand{\fin}{\mbox{\rm{Fin}}}
	\newcommand{\dom}{\mbox{\rm{dom}}}
	\newcommand{\ran}{\mbox{\rm{ran}}}
	\newcommand{\gwf}{g^{W,f}}
	\newcommand{\Sub}{\mbox{\rm{Sub}}}
	\newcommand{\Cl}{\mbox{{\rm Cl}}}
	\newcommand{\cont}{(2^{\omega})}
	\newcommand{\II}{\mbox{\rm{I}}}
	\newcommand{\Gr}{ $G=\langle V,E\rangle$ }
	\def\rel#1;{\mbox{rl}_G(#1)}
	\def\rls#1;{\mbox{rl}^*_G(#1)}
	\def\eqclass#1;#2;{#1/\equiv_{#1,#2}}
	\def\typ{\mbox{\rm tp${}_G$}}
	\def\tp{\typ}
	\def\eqr#1;{\equiv_{G,#1}}
	\def\noteqr#1;{\not\equiv_{G,#1}}
	\def\eqc#1;#2;{[#1]_{G,#2}}
	\def\oot{{{\omega}_2}}
	\def\twin#1;#2;{\mbox{twin${}_G$}(#1,#2)}
	\newcommand{\rank}{\mbox{rank}}
	\newcommand{\mod}{\mbox{ mod }}
	\def\pd#1;#2;{p_{{\delta}^{#1}_{#2}}}
	\def\gd#1;#2;{{\gamma}_{{\delta}^{#1}_{#2}}}
	\def\yyy#1;#2;#3;{{}^{#1}Y^{#2}_{#3}}
	\def\yy#1;#2;{{}^{#1}Y^{#2}}
	\def\dabb#1;#2;#3;#4;{\mbox{\rm D}^#1_#2(#3,#4)}
	\def\dab#1;#2;#3;{\mbox{\rm D}_#1(#2,#3)}
	\def\eps#1;{\epsilon_G(#1)}
	\def\ep#1;{\epsilon(#1)}
	\newcommand{\Prr}{{\em Pr}}
	\newcommand{\halfg}{[\oo;\oo]}
	\newcommand{\supp}{\mbox{\rm{supp}}}
	\newcommand{\nstat}{\mbox{\rm{NS}}}
	\newcommand{\height}{\mbox{\rm{height}}}
	\def\dt#1;{\mbox{\rm dt}_{#1}}
	\def\dtb{\dt b;}
	\newcommand{\compp}{\|_{{}_P}}

	\def\subs{\subset}
	\def\<{\left\langle}
	\def\>{\right\rangle}
	\def\cc{\gamma}
	\def\dd{\delta}
	\def\n{\nu}
	\def\m{\mu}
	\def\oo{\omega_1}
	\def\br#1;#2;{\bigl[ {#1} \bigr]^ {#2} }
	\def\bc#1;#2;{\bigl( {#1} \bigr)^ {#2} }
	\def\k{\kappa}
	\def\Dom{\operatorname{dom}}
	\def\Ran{\operatorname{ran}}
	\def\ooseq#1;#2;{\< {#1}_{#2}:{#2}<\oo\>}
	\def\ooset#1;#2;{\{ {#1}_{#2}:{#2}<\oo\}}
	\def\seq#1;#2;#3;{\< {#1}_{#2}:{#2}<#3\>}
	\def\set#1;#2;#3;{\{ {#1}_{#2}:{#2}<#3\}}
	\def\oseq#1;#2;{\< {#1}_{#2}:{#2}<\O\>}
	\def\oset#1;#2;{\{ {#1}_{#2}:{#2}<\O\}}
	\def\oosequ#1;#2;{\< {#1}^{#2}:{#2}<\oo\>}
	\def\oosetu#1;#2;{\{ {#1}^{#2}:{#2}<\oo\}}
	\def\sequ#1;#2;#3;{\< {#1}^{#2}:{#2}<#3\>}
	\def\setu#1;#2;#3;{\{ {#1}^{#2}:{#2}<#3\}}
	\def\osequ#1;#2;{\< {#1}^{#2}:{#2}<\O\>}
	\def\osetu#1;#2;{\{ {#1}^{#2}:{#2}<\O\}}
	\def\fsq#1;#2;{\<#1_0,\dots,#1_{#2}\>}
	\def\fsqv#1;#2;{\vec{#1}=\<#1_0,\dots,#1_{#2}\>}
	\def\con#1;{\operatorname{Con}(\text{#1})}
	\def\empt{\emptyset}
	\def\rest{\lceil}
	\def\dplus{\diamondsuit^{+}}

	\def\is{\cong}
	\def\bij{@>1-1>onto>}
	\def\abs#1{\|#1\|}
	\def\force{\raisebox{1.5pt}{\mbox{$\scriptscriptstyle\|$}}
	\mbox{$\!\mbox{---}$}}
	\def\qst{Q^{*}}

	\def\HH{\operatorname{Fin}}
	\def\inc{\not\|}

	\begin{abstract}
	Let $\kcal$ be the family of graphs on $\oo$ without cliques or independent
	subsets of size $\oo$. We prove that
	\begin{enumerate}\alabel
	\item it is consistent with CH that every $G\in\kcal$ has $2^{\oo}$ many
	pairwise non-isomorphic subgraphs,
	\item the following proposition holds in L:
	$(*)$ {\em there is a $G\in\kcal$ such that for each partition $(A,B)$ of $\oo$
	either $G\cong G[A]$ or  $G\cong G[B]$},
	\item the failure of $(*)$ is consistent with ZFC.
	\end{enumerate}
	\end{abstract}
	\section{Introduction}

	We assume only basic knowledge of set theory --- simple combinatorics
	for section \ref{sc:large}, believing in $L\models\dplus$ defined below
	for section
	\ref{sc:dplus}, and finite support iterated forcing for section
	\ref{sc:no_sm}.


	Answering a question of R. Jamison, H. A. Kierstead and P. J. Nyikos \cite{KN}
	proved that if  an $n$-uniform hypergraph $G=\<V,E\>$ is isomorphic
	to each of its induced subgraphs of cardinality $|V|$, then $G$ must be either
	empty or complete. They raised several new problems. Some of
them will be investigated in this paper. To present them we need to introduce some 
	notions.

	An infinite graph \Gr is called {\em non-trivial} iff $G$ contains no
	clique or independent subset of size $|V|$.
	Denote the class of all non-trivial graphs on $\oo$ by $\kcal$.
	Let $\II(G)$ be the set of all isomorphism classes
	of induced subgraphs of $G=\<V,E\>$ with size $|V|$.

	H. A. Kierstead and P. J. Nyikos proved that $|\II(G)|\ge{\omega}$ for
	each $G\in\kcal$  and asked whether $|\II(G)|\ge 2^{\omega}$
	or $|\II(G)|\ge 2^{\oo}$ hold or not.
	In \cite{HNS} it was shown that (i)  $|\II(G)|\ge 2^{\omega}$ for each
	$G\in\kcal$,
	(ii) under
	$\dplus$ there exists a $G\in\kcal$ with $|\II(G)|=\oo$.
	In section \ref{sc:large} we show that if ZFC is consistent, then so is
	ZFC + CH + ``{\em $|\II(G)|=2^{{\omega}_1}$ for each $G\in\kcal$}''.
	Given any $G\in\kcal$ we will investigate its
	partition tree. Applying  the weak $\diamondsuit$ principle of Devlin
	and Shelah \cite{DS} we show that if this partition tree is
	a special Aronszajn tree, then $|\II(G)|>\oo$.
	This result completes the investigation of problem 2 of \cite{KN}
	for ${{\omega}_1}$.

	Consider a graph \Gr . We say that $G$ is {\em almost smooth} if it is
	isomorphic to $G[W]$ whenever $W\subset V$ with $|V\setminus W|<|V|$.
	The graph $G$ is called {\em quasi smooth} iff it is isomorphic either to
	$G[W]$ or to $G[{V\setminus W}]$ whenever $W\subset V$.
	H. A. Kierstead and P. J. Nyikos asked (problem 3) whether an
	almost smooth,
	non-trivial graph can exist.
	In \cite{HNS} various models of ZFC  was constructed which contain
	such graphs on ${{\omega}_1}$. It was also shown
	that the existence of a non-trivial, quasi smooth graph on
	${{\omega}_1}$ is consistent with ZFC. But in that model CH failed.
	In section \ref{sc:dplus}
	we prove that $\dplus$, and so V=L, too,  implies the existence of such
	a graph.

	In  section \ref{sc:no_sm} we construct a model of ZFC
	in which there is no quasi-smooth  $G\in\kcal $.
	Our main idea is that given a $G\in\kcal$ we try to  construct
	a  partition $(A_0,A_1)$ of $\oo$ which is so bad that not only
	$G\not\cong G[A_i]$ in the ground model
	 but certain simple generic extensions   can not add such
isomorphisms to the ground model. 
	We  divide the class $\kcal$ into three subclasses
	and develop different methods to carry out our plan.

	The question whether the existence of an
	almost-smooth  $G\in\kcal$
	can be proved in ZFC is still open.

	We  use the standard set-theoretical notation throughout, cf \cite{J}.
	Given a graph \Gr we write  $V(G)=V$ and $E(G)=E$. If
	 $H\subset V(G)$ we define $G[H]$ to be
	$\langle H,E(G) \cap [H]^2\rangle $. Given $x\in V $
	take $G(x)=\bigl\{y\in V:\{x,y\}\in E\bigr\}$.
	If $G$ and $H$ are graphs we write $G\cong H$ to mean that $G$ and $H$
	are isomorphic. If $f:V(G)\to V(H)$ is a function we denote by $f:G\cong H$ the
	fact
	that $f$ is an isomorphism between $G$ and $H$.

	Given a set $X$ let $\bijp(X)$ be the set of all bijections between subsets of
	$X$. If \Gr is a graph take
	\[
	\isop(G)=\left\{f\in\bijp(V):\ f:G[\dom(f)]\cong G[\ran(f)]\right\}.
	\]
	We denote by $\fin(X,Y)$ the set of all functions mapping a finite subset of
	$X$ to $Y$.

	Given a poset $P$ and $p,q\in P$ we write $p\compp q$ to mean that
	$p$ and $q$ are compatible in $P$.

	The axiom  $\diamondsuit^{+}$ claims that {\em
	there is a sequence
	$\bigl\langle
	 S_{\alpha}:{\alpha}<{{\omega}_1}\bigr\rangle$ of contable sets
	such that for each $X\subset {{\omega}_1}$ we have a closed unbounded
	$C\subset{{\omega}_1}$ satisfying $X\cap{\nu}\in S_{\nu}$ and
	$C\cap\nu\in S_{\nu}$ for each $\nu\in C$.}

	We denote by $\mbox{{\rm TC}}(x)$ the transitive closure of a set $x$.
	  If ${\kappa}$ is a cardinal take
	$H_{\kappa}=\{x:|\mbox{{\rm TC}}(x)|<{\kappa}\}$ and
	$\hcal_{{\kappa}}=\langle H_{\kappa},\in\rangle$.

	Let us denote by  $\dcal_{\oo}$ the club filter on $\oo$.

	\section{ I($G$) can be always large}
	\label{sc:large}
	\begin{theorem}
	\label{th:ar}
	Asume that  GCH holds and every Aronszajn-tree is special. Then
	$|\II(G)|=2^{{\omega}_1}$ for each $G\in\kcal$.
	\end{theorem}

	\begin{demo}{Remark} S.Shelah proved,
	\cite[chapter V. \S  6,7]{S}, that the assumption
	of theorem \ref{th:ar} is consistent with ZFC.
	\end{demo}
	During the proof we will apply the following definitions and lemmas.
	\begin{lemma} Assume that $G\in\kcal$,
	$A\in[{{\omega}_1}]^{{\omega}_1}$ and
	$|\{G(x)\cap A:x\in{{\omega}_1}|={{\omega}_1}$. Then $|\II(G)|=2^{{\omega}_1}$.
	\end{lemma}

	\proof See \cite[theorem 2.1 and lemma 2.13]{HNS}.
	\eproof

	\begin{definition}
	\label{df:ptree}
	Consider a graph $G=\left\langle {{\omega}_1},E\right\rangle $.
	\begin{enumerate}
	\item For each ${\nu}\in{{\omega}_1}$ let us define the ordinal
	${\gamma}_{\nu}\in{{\omega}_1}$ and the sequence
	$\left\langle {\xi}^{\nu}_{\gamma}:{\gamma}\le{\gamma}_{\nu}\right\rangle $
	as follows:  put $\xi^{\nu}_0=0$ and if
	$\left\langle {\xi}^{\nu}_{\alpha}:{\alpha}<{\gamma}\right\rangle $
	is defined, then take
	\begin{displaymath}
	{\xi}^{\nu}_{\gamma}=
	     \min\left\{{\xi}:\forall {\alpha}<{\gamma}\ {\xi}>{\xi}^{\nu}_{\alpha}
	         \mbox{ and }
	              (\{{\xi}^{\nu}_{\alpha},{\xi}\}\in E\mbox{ iff }
	              \{{\xi}^{\nu}_{\alpha},{\nu}\}\in E)
	           \right\}.
	\end{displaymath}
	If ${\xi}^{\nu}_{\gamma}={\nu}$, then we put ${\gamma}_{\nu}={\gamma}$.
	\item Given ${\nu},{\mu}\in{{\omega}_1}$ write
	${\nu}\prec^G{\mu}$ iff
	${\xi}^{\nu}_{\gamma}={\xi}^{\mu}_{\gamma}$
	for each ${\gamma}\le{\gamma}_{\nu}$.
	\item Take $\tcal^G=\left\langle {{\omega}_1},\prec^G \right\rangle $.
	$\tcal_G$ is called the {\em partition tree of $G$}.
	\end{enumerate}
	\end{definition}

	\begin{lemma}
	If $G=\left\langle {{\omega}_1},E\right\rangle \in\kcal$
	with $|\II(G)|<2^{{\omega}_1}$, then $\tcal^G$ is an Aronszajn tree.
	\end{lemma}

	\proof
	By the construction of $\tcal^G$, if ${\nu},{\mu}\in{{\omega}_1}$,
	${\nu}<{\mu}$ and $G({\nu})\cap{\nu}=G({\mu})\cap {\nu}$, then
	${\nu}\prec^G{\mu}$. So the levels of $\tcal^G$ are countable by lemma 2.2.
	On the other hand, $\tcal^G$
	does not contain ${{\omega}_1}$-branches, because the branches are
	prehomogeneous subsets and $G$ is non-trivial.
	\eproof
	\begin{definition}

	\begin{enumerate}
	\item Let $F:\cont^{<{{\omega}_1}}\to 2$ and $A\subset{{\omega}_1}$.
	We say that a function $g:{{\omega}_1}\to 2$ is an {\em $A$-diamond for $F$}
	iff, for any $h\in\cont^{{\omega}_1}$,
	$\left\{{\alpha}\in A: F(h\rest{\alpha})=g({\alpha})\right\}$
	is a stationary subset of ${{\omega}_1}$.
	\item $A\subset{{\omega}_1}$ is called {\em a small subset of ${{\omega}_1}$}
	iff for some $F:\cont^{<{{\omega}_1}}\to 2$ no function is an
	$A$-diamond for $F$.
	\item $\jcal=\left\{A\subset{{\omega}_1}:A
	\mbox{ is a small subset of ${{\omega}_1}$}
	\right\}$.
	\end{enumerate}
	\end{definition}
	In \cite{DS} the following was proved:
	\begin{theorem}
	\label{th:ds}
	If $2^{\omega}<2^{{\omega}_1}$, then $\jcal$ is a countably complete, proper,
	normal ideal on ${{\omega}_1}$.
	\end{theorem}

	After this preparation  we are ready to prove theorem
	\ref{th:ar}.
	\proof
	Assume that $G=\left\langle {{\omega}_1},E\right\rangle\in\kcal $.

	$|\II(G)|<2^{{\omega}_1}$ and a contradiction will be derived.

	Since $2^{{\omega}_1}={\omega}_2$, we can fix a sequence
	$\left\{G_{\nu}:{\nu}<{{\omega}_1}\right\}$ of graphs on ${{\omega}_1}$
	such that for each $Y\in[{{\omega}_1}]^{{\omega}_1}$ there is a
	${\nu}<{{\omega}_1}$ with $G[Y]\cong G_{\nu}$.
	Write $G_{\nu}=\langle \omega_1,E_{\nu}\rangle$.

	\vspace{3pt}
	Consider the Aronszajn-tree
	$\tcal^G=\left\langle {{\omega}_1},\prec^G \right\rangle $.
	Since every Aronszajn-tree is special and $\ical$ is a countably
	complete ideal on ${{\omega}_1}$, there is an antichain $S$ in $\tcal^G$
	with $S\notin\jcal$.
	Take
	\[
	A=\left\{{\alpha}\in{{\omega}_1}
	   :\exists{\sigma}\in S ({\alpha}\prec^G {\sigma})\right\}.
	\]
	Now property $(*)$ below holds:
	\[
	\makebox[2.5cm][l]{$(*)$}
	\forall {\sigma}\in S\ \forall{\rho}\in (S\cup A)\setminus{{\sigma}+1}
	\makebox[5cm]{}
	\]
	\vspace{-0.7cm}
	\[
	\makebox[4cm]{}
	\exists{\alpha}\in A\cap{\sigma}\ (\{{\sigma},{\alpha}\}\in E
	\mbox{ iff } \{{\rho},{\alpha}\}\notin E).
	\]

	Indeed, if for each ${\alpha}\in A\cap{\sigma}$ we had
	$\{{\sigma},{\alpha}\}\in E$
	 iff $\{{\rho},{\alpha}\}\in E$, then $\sigma\prec^{G}{\rho}$ would hold
	 by the construction of $\tcal^G$.

	Let ${\nu}\in{{\omega}_1}$, ${\sigma}\in S$, $T\subset S\cap {\sigma}$
	and $f:G[(A\cap{\sigma})\cup T]\to G_{\nu}$ be an embedding.
	Define $F({\nu},{\sigma},T,f)\in 2$ as follows:
	\[
	F({\nu},{\sigma},T,f)=1 \mbox{ iff }\exists x\in G_{\nu}
	(\forall{\alpha}\in A\cap{\sigma})
	\bigl(\{x,f({\alpha})\}\in E_{\nu}\mbox{ iff }
	\{{\sigma},{\alpha}\}\in E\bigr).
	\]

	In case ${\omega}{\sigma}={\sigma}$, under suitable encoding,
	$F$ can be viewed as a function from $\cont^{<{\omega}_1}$
	to $2$.

	Since $S\notin\jcal$, there is a $g\in 2^{{\omega}_1}$ such that
	for every ${\nu}\in{{\omega}_1}=2^{\omega}$, $T\subset S$
	and $f:G[A\cup T]\cong G_{\nu}$, the set
	\[
	S_T=\left\{{\sigma}\in S:g({\sigma})=
	  F({\nu},{\sigma},T\cap{\sigma},f\rest{\sigma})\right\}
	\]
	is stationary.
	Take $T=\{{\sigma}\in S:g({\sigma})=0\}$.
	Choose an ordinal ${\nu}<{{\omega}_1}$ and a function $f$
	with $f:G[A\cup T]\cong G_{\nu}$.
	For each
	${\sigma}<{{\omega}_1}$ with ${\sigma}={\omega}{\sigma}$ it follows,
	by $(*)$,  that
	\[
	{\sigma}\in T \mbox{ iff  }\exists x\in \omega_1
	\ \forall{\alpha}\in S\cap{\sigma}
	\ (\{x,f({\alpha})\}\in E_{\nu}\mbox{ iff }\{{\sigma},{\alpha}\}\in E).
	\]
	 Thus $g({\sigma})=0$ iff
	$F({\nu},{\sigma},T\cap {\sigma},f\rest {\sigma})=1$,
	for each $\sigma\in S$, that is,
	$S_T=\emptyset$, which is a contradiction.
	\eproof

	\section{A quasi-smooth graph under $\dplus$}
	\label{sc:dplus}
	\begin{theorem}
	\label{th:dplus}
	If $\dplus$ holds, then there exists a non-trivial, quasi-smooth graph on $\oo$.
	\end{theorem}

	\proof
	Given a set $X$, $\acal{\subset} \mbox{P}(X)$ and $\fcal{\subset}\bijp(X)$ take

	\[
	\Cl(\acal,\fcal)=\bigcap\left\{\bcal:\bcal\supset\acal\mbox{ and }
	\forall B_0,B_1\in\bcal\ \forall f\in\fcal\ \forall Y\in [X]^{<{\omega}}
	\right.\makebox[1cm]{}
	\]
	\vspace{-0.8cm}
	\[
	\makebox[6cm]{} \left.\left\{ B_0\cup B_1,f''B_0,B_0\triangle
Y \right\} {\subset}\bcal \right\}. 
	\]
	We say that $\acal$ is {\em $\fcal$-closed} if $\acal=\Cl(\acal,\fcal)$.
	Given $\acal,\dcal{\subset} P(X)$, we say that $\dcal$ is {\em
uncovered by $\acal$} 
	if $|D{\setminus} A|={\omega}$ for each $A\in\acal$ and $D\in\dcal$.

	\begin{lemma}
	\label{lm:partition}
	Assume that $\fcal{\subset}\bijp(X)$ is a countable set,
	$\acal^0$, $\acal^1{\subset} \mbox{P}(X)$ are countable,
	$\fcal$-closed families.
	If $\dcal{\subset} P(X)$ is a countable family which is uncovered by
	$\acal^0\cup\acal^1$, then there is a partition $(B_0,B_1)$ of $X$
	such that $\dcal$ is uncovered by
	$\Cl(\acal^i\cup \left\{ B_i \right\} ,\fcal)$
	for $i<2$.
	\end{lemma}

	\proof
	We can assume that $\fcal$ is closed under composition. Fix an enumeration
	$ \left\{ \left\langle D_n,k_n,F_n,i_n,A_n\right\rangle :n\in{\omega} \right\} $
	 of
	$\dcal\times{\omega}\times\fcal^\fini\times
	\left\{ \left\langle i,A\right\rangle :i\in 2, A\in \acal^i \right\} $.
	By induction on $n$, we will pick points $x_n\in X$ and will define finite sets,
	$B^0_n$ and $B^1_n$, such that $B^0_n\cap B^1_n={\emptyset}$ and
	$B^i_n{\subset} B^i_{n+1}$.

	Assume that we have done it for $n-1$. Write
	$F_n=\left\langle f_0,\ldots,f_{k-1}
	\right\rangle $. Take $B_{n-1}=B^0_{n-1}\cup B^1_{n-1}$
	and
	\begin{displaymath}
	B^-_n=B_{n-1}\cup\bigcup \left\{ f_j''B_{n-1}:j<k \right\}.
	\end{displaymath}
	Pick an arbitrary point $x_n\in D_n{\setminus}(A_n\cup B^-_n)$.
	Put
	\begin{displaymath}
	B^{i_n}_n=B^{i_n}_{n-1}
	\end{displaymath}
	and
	\begin{displaymath}
	B^{1-i_n}_n=B^{1-i_n}_{n-1}\cup \left\{ x_n \right\} \cup \left\{ f_j^{-1}(x_n):
	j<k \right\} .
	\end{displaymath}
	Next choose a partition $(B^0,B^1)$ of $X$ with
	$B^i\supset\cup \{ B^i_n:n<{\omega} \} $ for $i<2$.
	We claim that it works.
	Indeed, a typical element of
	$\Cl(\acal^i\cup \{ B^i \} ,\fcal)$ has the form
	\begin{displaymath}
	C=A\cup\bigcup \left\{ f_j''B^i:j<k \right\} ,
	\end{displaymath}
	where $A\in\acal$, $k<{\omega}$ and $f_0,\ldots,f_{k-1}\in\fcal$.
	So, if $D\in\dcal$, then
	\begin{displaymath}
	D{\setminus} C\supset \left\{ x_n:D_n=D,A_n=A, i_n=i \mbox{
and } F_n=\left\langle f_0,\ldots,f_{k-1}\right\rangle \right\} 
	\end{displaymath}
	because $x_n\notin A$ and $f_j^{-1}(x_n)\in B^{1-i}$ by the constuction.
	\eproof
	Consider a sequence $F=\left\langle f_0,\ldots,f_{n-1} \right\rangle $.
	 Given a family $\fcal{\subset}\bijp(X)$ we say that $F$ is an {\em
	$\fcal$-term}
	provided $f_i=f$ or $f_i=f^{-1}$ for some $f\in\fcal$, for each $i<n$.
	We denote the function $f_0\circ\cdots\circ f_{n-1}$ by $F$ as well.
	We will assume that the empty term denotes the identity function on $X$.
	If $l\le n$ take $_{(l)}F=\left\langle f_0,\ldots,f_{l-1} \right\rangle $ and
	$F_{(l)}=\left\langle f_l,\ldots,f_{n-1} \right\rangle $.
	Let
	\begin{displaymath}
	\Sub(F)= \left\{ \left\langle f_{i_0},\ldots,f_{i_{l-1}} \right\rangle :
	l\le n, i_0<\ldots< i_{l-1}<n
	 \right\} .
	\end{displaymath}
	Given $f\in\fcal$ and $x,y\in X$ with $x\notin \dom(f)$
	and $y\notin \ran(f)$ let $F^{f,x,y}$ be the term
	that we obtain replacing each occurrence of $f$ and of $f^{-1}$ in $F$ with
	$f\cup \left\{ \left\langle x,y \right\rangle  \right\} $
	and with $f^{-1}\cup \left\{ \left\langle y,x \right\rangle  \right\} $, respectively.
	\begin{lemma}
	\label{lm:fxy}
	Assume that $\fcal{\subset}\bijp(X)$, $\acal{\subset} P(X)$ is $\fcal$-closed,
	$F_0,\ldots,F_{n-1}$ are $\fcal$-terms, $z_0,\dots,z_{n-1}\in X$,
	$A_0,\ldots, A_{n-1}\in\acal$
	such that
	for each $i<n$
	\begin{displaymath}
	\makebox[4cm][l]{$(*)$} z_i\notin\bigcup \left\{ F''A_i:F\in\Sub(F_i) \right\}.
	\makebox[4cm][l]{}
	\end{displaymath}
	If $f\in\fcal$, $x\in X{\setminus}\dom(f)$, $Y\in[X{\setminus}\ran(f)]^{\omega}$
	with $|A\cap Y|<{\omega}$
	for each $A\in\acal$, then there are infinitely many $y\in Y$
	such that $(*)$ remains true when
	replacing $f$ with $f\cup \left\{  \left\langle
	x,y \right\rangle  \right\} $,
	that is,
	\begin{displaymath}
	\makebox[3.8cm][l]{$(**)$}
	z_i\notin\bigcup \left\{ F''A_i:F\in\Sub(F^{f,x,y}_i) \right\}
	\makebox[3.8cm][l]{}
	\end{displaymath}
	for each $i<n$.
	\end{lemma}

	\proof
	It is enough to prove it for $n=1$.
	Write $F= \left\langle f_0\ldots,f_{k-1}\right\rangle$, $A=A_0$, $z=z_0$.
	Take
	\begin{displaymath}
	Y_{F,A}= \left\{ y\in Y: \mbox{$(**)$ holds for $y$} \right\} .
	\end{displaymath}

	Now we prove the lemma by induction on $k$.

	If $k=0$, then
	$Y_{F,A}=Y{\setminus} A$.
	Suppose we know the lemma for $k-1$.   Using the induction hypothesis
	we can assume that
	(\dag) below holds:
	\begin{displaymath}
	\makebox[2cm][l]{(\dag)}
	Y=\bigcap \left\{ Y_{G,F_{(l)}''A}:l\le n,G\in \Sub({}_{(l)}F^{f,x,y}),
	G\ne F^{f,x,y} \right\}.
	\makebox[2cm][l]{}
	\end{displaymath}

	Assume  that $|Y_{F,A}|<{\omega}$
	and a contradiction will be derived.

	First let us remark that either $f_{k-1}=f$ or $f_{k-1}=f^{-1}$
	by (\dag).

	\begin{case}{1}{$f_{k-1}=f^{-1}$.}
	Then $Y_{F,A}{\supset} Y{\setminus} A$ by (\dag), so we are done.
	\end{case}
	\begin{case}{2}{$f_{k-1}=f$.}
	In this case  $x\in A$ and for all but finitely many $y\in Y$
	we have $z=F^{f,x,y}(x)$. Then for each $y,y'\in Y$ take
	\begin{displaymath}
	l(y,y')=\max \left\{ l\le n:
	\forall i<l\ F_{(i)}^{f,x,y}(x)=F_{(i)}^{f,x,y'}(x) \right\} .
	\end{displaymath}
	By Ramsey's theorem, we can assume that $l(y,y')=l$ whenever $y,y'\in Y$.
	Clearly $l<n$. Then
	$F_{(l)}^{f,x,y}(x)\ne F_{(l)}^{f,x,y'}(x)$ but
	$F_{(l-1)}^{f,x,y}(x)=F_{(l-1)}^{f,x,y'}(x)$, so $f_l=f^{-1}$ and
	$F_{(l-1)}^{f,x,y}(x)=x$ for each $y\in Y$.
	Thus $z={}_{(l-1)}F^{f,x,y}(x)$ for each $y\in Y$, which contradicts
	(\dag) because $x\in A$.
	\end{case}
	The lemma is proved.
	\eproof

	We are ready to construct our desired graph.

	First fix a sequence $ \left\langle M_{\alpha}:{\alpha}<\oo \right\rangle $ of countable, elementary
	submodels of some $\hcal_{\lambda}$  with
	$\left\langle M_{\gamma}:{\gamma}<{\alpha} \right\rangle \in M_{\alpha}$
	for each ${\alpha}<\oo$,
	where ${\lambda}$ is a large enough
	regular cardinal.

	Then choose a
	$\diamondsuit$-sequence
	$ \left\langle S_{\alpha}:{\alpha}<\oo \right\rangle \in M_0$
	for the uncountable subsets of $\oo$,
	that is ,
	$ \left\{ {\alpha}<{{\omega}_1}:X\cap{\alpha}=S_{\alpha} \right\}
	  \notin NS({\omega_1})$
	whenever $X\in[\oo]^{\oo}$. We can also assume that $S_{\alpha}$
	is cofinal in ${\alpha}$ for each limit ${\alpha}$.

	We will define, by induction on ${\alpha}$,
	\begin{enumerate}
	\item graphs $G_{\alpha}= \left\langle {\omega}{\alpha},E_{\alpha} \right\rangle	 $
	with $G_{\beta}=G_{\alpha}[{\omega}{\beta}]$ for ${\beta}<{\alpha}$,
	\item countable sets $\fcal_{\alpha}\in\isop(G_{\alpha})$,
	\end{enumerate}
	satisfying the induction hypotheses (I)--(II) below:
	\begin{enumerate} \Rlabel
	\item $ \left\{ S_{{\omega}{\gamma}}:{\gamma}\le{\alpha} \right\} $ is uncovered
	 by
	$I_{\alpha}\cup J_{\alpha}$
	where
	\begin{displaymath}
	I_{\alpha}=
	   \Cl( \left\{ G({\nu})\cap{\nu}:{\nu}\in{\omega}{\alpha} \right\} ,
	   \bigcup\limits_{{\beta}\le{\alpha}}\fcal_{\beta})
	\end{displaymath}
	and
	\begin{displaymath}
	J_{\alpha}=
	   \Cl( \left\{ {\nu}{\setminus} G({\nu}):{\nu}\in{\omega}{\alpha} \right\} ,
	   \bigcup\limits_{{\beta}\le{\alpha}}\fcal_{\beta}).
	\end{displaymath}
	\end{enumerate}
	To formulate (II)
	we need the following definition.
	\begin{definition}
	\label{df:large}
	Assume that ${\alpha}={\beta}+1$ and $Y{\subset}{\omega}{\alpha}$.
	We say that $Y$ is {\em large} if
	$\forall n\in{\omega}$, $\forall \left\langle  \left\langle f_i,x_i \right\rangle :i<n \right\rangle $,
	$\forall h$ \newline {\bf if}
	\begin{enumerate}
	\item $\forall i<n\ \exists {\alpha}_i<{\beta}\ f_i\in\fcal_{{\alpha}_i}$,
	\item $\forall i<n\ {\omega}{\alpha}_i\le x_i<{\omega}{\beta}$,
	\item $\forall i<n\ \ran(f_i){\subset} Y$,
	\item $\forall i\ne j<n\ \ran(f_i)\cap \ran(f_j)=\emptyset$
	\item $h\in \mbox{{\rm Fin}}(Y\cap {\omega}{\beta},2)$ and
	$\dom(h)\cap\bigcup \left\{ \ran(f_i):i<n \right\} =\emptyset$,
	\end{enumerate}
	{\bf then}
	\newline
	$\exists y\in Y\cap[{\omega}{\beta},{\omega}{\alpha})$
	such that
	\begin{enumerate}\addtocounter{enumi}{5}
	\item $\forall i<n\ \forall x\in\dom(f_i)
	\ (\left\{ y,f_i(x) \right\} \in E_{\alpha}
	\mbox{ iff }  \left\{ x_i,x \right\} \in E_{\alpha})$,
	\item $\forall z\in\dom(h)\ \left\{ y,z \right\} \in E_{\alpha}\mbox{ iff }
	h(z)=1$.
	\end{enumerate}
	\end{definition}

	Take
	\begin{enumerate}\addtocounter{enumi}{1} \Rlabel
	\item If ${\alpha}={\beta}+1$, then ${\omega}{\alpha}$ is large.
	\end{enumerate}
	The construction will be carried out in such a way that
	$ \left\langle G_{\beta}:{\beta}\le{\alpha} \right\rangle \in M_{\alpha}$
	and $ \left\langle \fcal_{\beta}:{\beta}<{\alpha} \right\rangle \in M_{\alpha}$.

	To start with take $G_0= \left\langle \emptyset,\emptyset
\right\rangle $ and $\fcal= \left\{ \emptyset \right\} $. 
	Assume that the construction is done for ${\beta}<{\alpha}$.

	\begin{case}{1}{${\alpha}$ is limit.}
	We must take $G_{\alpha}=\cup \left\{ G_{\beta}:{\beta}<{\alpha} \right\} $.
	We will define sets
	$\fcal^0_{\alpha}, \fcal^1_{\alpha}{\subset}\isop(G_{\alpha})$
	and will take $\fcal_{\alpha}=\fcal^0_{\alpha}\cup\fcal_{\alpha}^1$.

	Let
	\[
	         \fcal^0_{\alpha}=
	        \left\{ f\in\isop(G_{\alpha})\cap M_{\alpha}:
	          \exists \left\langle {\alpha}_n
	:n<{\omega}\right\rangle
	              {\subset}{\alpha}
	   \ \sup \left\{ {\alpha}_n:n<{\omega} \right\} ={\alpha},
	        \right.
	\makebox[1cm]{}
	\]
	\vspace{-0.7cm}
	\[
	\makebox[2cm]{}
	      \left.f\rest{\omega}{\alpha}_n\in\fcal_{{\alpha}_n}
	             \mbox{ and }
	       f\rest{\omega}{\alpha}_n:G_{{\alpha}_n}\cong G_{{\alpha}_n}[\ran(f)]
	       \mbox{ for each $n\in{\omega}$} \right\}.
	\]
	Take
	$
	\fcal^-=\bigcup\limits_{{\beta}<{\alpha}}\fcal^{\beta}\cup\fcal^0_{\alpha}$,
	$I^-_{\alpha}=\bigcup\limits_{{\beta}<{\alpha}}I_{\beta}$ and
	$J^-_{\alpha}=\bigcup\limits_{{\beta}<{\alpha}}J_{\beta}$.
	Clearly $\fcal^-{\subset} M_{\alpha}$ with $\fcal^-\in M_{{\alpha}+1}$, so
	$M_{{\alpha}+1}\models\mbox{``}|\fcal^-|={\omega}\mbox{''}$.
	Obviously both $I^-_{\alpha}$ and $J^-_{\alpha}$ are $\fcal^-$-closed
	and $\scal= \left\{ S_{{\omega}{\beta}}:{\beta}\le{\alpha}
\right\} $ is uncovered by 
	them.

	>From now on we work in $M_{{\alpha}+1}$ to construct $\fcal^1_{\alpha}$.
	For $W{\subset}{\omega}{\alpha}$ write
	$L_W= \{ {\nu}<{\alpha}:W\cap({\omega}{\nu}+{\omega})\mbox{ is large}\}.$

	Take
	\[
	\wcal_{\alpha}=
	    \left\{
	    \left\langle W,f \right\rangle
	         \in(P({\omega}{\alpha})\cap\mcal_{\alpha})
	\times
	                   (\cup_{{\beta}<{\alpha}}\fcal_{\beta}):
	             \mbox{$L_W$ is cofinal in ${\alpha}$}\right.
	\makebox[1cm]{}
	\]
	\vspace{-0.6cm}
	\[
	\makebox[2cm]{}
	     \left. \makebox[3cm][r]{and }
	f:G_{{\gamma}_f}\cong G_{{\gamma}_f}[W\cap {\omega}{\gamma}_f]
	      \mbox{ for some ${\gamma}_f<{\alpha}$}
	  \right\}.
	\]

	We want to  find functions $g^{W,f}\supset f$
	for $ \left\langle W,f \right\rangle \in\wcal_{\alpha}$ such that
	\begin{enumerate}\Alabel
	\item $g^{W,f}:G_{\alpha}\cong G_{\alpha}[W]$
	\item taking $\fcal^1_{\alpha}= \left\{ g^{W,f}: \left\langle W,f \right\rangle
	\in\wcal_{\alpha} \right\} $
	the induction hypothesis (I) remains true.
	\end{enumerate}

	First we prove a lemma:

	\begin{lemma}
	\label{lm:dense}
	If $ \left\langle W,f \right\rangle \in\wcal_{\alpha}$, $g\in\isop(G_{\alpha}, G_{\alpha}[W])$,
	$g\supset f$, $|g{\setminus} f|<{\omega}$, then
	\begin{enumerate}\rlabel
	\item for each $x\in W{\setminus}\dom(f)$ the set
	\[
	\left\{ y\in W:g\cup \left\{  \left\langle x,y \right\rangle  \right\} \in\isop(
	G_{\alpha},G_{\alpha}[W]) \right\}
	\]
	is cofinal in ${\omega}{\alpha}$.
	\item for each $y\in W{\setminus}\ran(f)$ the set
	\[
	 \left\{ x\in W:g\cup \left\{  \left\langle x,y \right\rangle  \right\} \in\isop
	(G_{\alpha},G_{\alpha}[W]) \right\}
	\]
	is cofinal in ${\omega}{\alpha}$.
	\end{enumerate}
	\end{lemma}

	\proof
	(i): Define the function $h:\ran(g){\setminus} ran(f)\to 2 $
	with $h(g(z))=1$ iff $ \left\{ z,x \right\} \in E_{\alpha}$. Choose ${\beta}\in
	L_W$
	with $\ran(h){\subset}{\omega}{\beta}$ and ${\gamma}_f\le{\beta}$. Since
	$W\cap({\omega}{\beta}+{\omega})$ is large, we have a
	$y\in W\cap[{\omega}{\beta},{\omega}{\beta}+{\omega})$ such that
	\begin{enumerate}
	\item  $ \left\{ y,f(z) \right\} \in E_{\alpha}$ iff
	$ \left\{ x,z \right\} \in E_{\alpha}$   for each
	$z\in\dom(f)$
	\item $ \left\{ y,g(z) \right\} \in E_{\alpha}$ iff
$h(g(z))=i$ for each $z\in\dom(g){\setminus}\dom(f)$. 
	\end{enumerate}
	But this means that $g\cup \left\{  \left\langle x,y \right\rangle  \right\} \in
	\isop(G_{\alpha},G_{\alpha}[W])$.
	\newline
	(ii) The same proof works using that ${\omega}{\beta}+{\omega}$ is large
	for each ${\beta}<{\alpha}$.
	\eproof

	By induction on $n$, we will pick points $z_n\in{\omega}{\alpha}$
	and will construct families of partial automorphisms,
	$ \left\{ \gwf_n: \left\langle W,f \right\rangle
\in\wcal_{\alpha} \right\} $ such that 
	$\gwf=\cup \left\{ \gwf_n:n<{\omega} \right\} $ will work.

	During the inductive construction we will speak about
	$\fcal_{\alpha}$-terms and about functions which are represented by them
	in the $n^{\mbox{\tiny th}}$ step.

	If $F= \left\langle h_0,\ldots h_{k-1} \right\rangle $ is an
	$\fcal_{\alpha}$-term and $n\in{\omega}$
	take $F_{[n]}=j_0\circ\cdots\circ j_{k-1}$ where
	\begin{displaymath}
	j_i= \left\{
	\begin{array}{ll}
	\gwf_n&\mbox{if $h_i=\gwf$,}\\
	(\gwf_n)^{-1}&\mbox{if $h_i=(\gwf)^{-1}$,}\\
	h_i&\mbox{otherwise.}
	\end{array}
	\right.
	\end{displaymath}
	First fix an enumeration
	$ \left\{  \left\langle  \left\langle W_n,f_n \right\rangle ,u_n,i_n
	\right\rangle :1\le n<{\omega} \right\} $ of
	$\wcal_{\alpha}\times{\omega}{\alpha}\times 2$
	and an enumeration
	$\left\langle\left\langle  \left\langle F_{n,i}:i<l_n \right\rangle
	,j_n, \left\langle A_{n,i}:i<l_n \right\rangle , D_n \right\rangle
	n<{\omega}\right\rangle$
	of
	the  quadruples
	$ \left\langle  \left\langle F_0,\ldots,F_{k-1} \right\rangle ,j, \left\langle
	A_0,\ldots,A_{k-1} \right\rangle , D \right\rangle $
	where $k<{\omega}$, $F_0,\ldots,F_{k-1}$ are $\fcal_{\alpha}$-terms,
	$j\in 2$, $D\in\scal$ and either $j=0$ and $A_0,\ldots,A_{k-1}\in I^-_{\alpha}$
	or
	$j=1$ and $A_0,\ldots,A_{k-1}\in J^-_{\alpha}$.

	During the inductive construction conditions (i)--(v) below will be satisfied:
	\begin{enumerate}\rlabel
	\item $\gwf_0=f$
	\item $\gwf_n\in\isop(G_{\alpha},G_{\alpha}[W])$
	\item $\gwf_n\supset\gwf_{n-1}$, $|\gwf_n{\setminus} f|<{\omega}$
	\item $z_k\notin\bigcup \left\{ F_{[n]}''A_{k,i}:F\in\Sub(F_{k,i}) \right\} $
	for each $i<l_k$ and $k<n$
	\item if $i_n=1$, then $u_n\in\dom(g^{W_n,f_n}_n)$,
	\newline
	if $i_n=0$, then either $u_n\notin W_n$ or $u_n\in\ran(g^{W_n,f_n}_n)$.
	\end{enumerate}

	If $n=0$, then take $\gwf_0=f$.

	If $n>0$, then let $\gwf_n=\gwf_{n-1}$
	whenever $ \left\langle W,f \right\rangle \ne \left\langle W_n,f_n \right\rangle	 $.
	Assume that $i_n=0$,
	$ \left\langle W,f \right\rangle =\left\langle W_n,f_n
\right\rangle $ and $u_n\notin\dom(g^{W_n,f_n}_{n-1})$. Then, by lemma
\ref{lm:dense},  
	 the set
	$Y= \left\{ y\in W:\gwf_{n-1}\cup \left\{  \left\langle u_n,y
\right\rangle  \right\} \in\isop(G_{\alpha},G_{\alpha}[W]) \right\} $ 
	is unbounded in ${\omega}{\alpha}$.
	Since  the members of $I^-_{\alpha}\cup J^-_{\alpha}$
	are bounded in ${\omega}{\alpha}$, we can apply lemma \ref{lm:fxy}
	to pick a point $y\in Y$ such that taking
	$g^{W_n,f_n}_n=g^{W_n,f_n}_{n-1}\cup \left\{  \left\langle u_n,y \right\rangle
	\right\} $  condition (iv)
	holds.

	If $i_n=1$ and
	$ \left\langle W,f \right\rangle =\left\langle W_n,f_n \right\rangle $,
	then the same argument works.

	Finally pick a point
	\begin{displaymath}
	z_n\notin D_n{\setminus}\bigcup \left\{ F_{[n]}''A_{n,i}:
	F\in\Sub(F_{n,i}),i<l_n\right\}.
	\end{displaymath}
	The inductive construction is done.

	Take $\gwf=\cup \left\{ \gwf_n:n<{\omega} \right\} $. By (v),
	\begin{displaymath}
	\gwf_n:G_{{\omega}{\alpha}}\cong G_{{\omega}{\alpha}}[W].
	\end{displaymath}
	By (iv), we have
	\begin{displaymath}
	z_k\in D_k{\setminus}\bigcup \left\{ F_{k,i}''A_{k,i}:i<l_k \right\}
	\end{displaymath}
	and so it follows that $ \left\{ S_{{\omega}{\beta}}:{\beta}\le{\alpha} \right\}
	 $
	is uncovered by $I_{\alpha}\cup J_{\alpha}$.

	\end{case}
	\begin{case}{2}{${\alpha}={\beta}+1$.}

	To start with we fix an enumeration
	$ \left\{
	\left\langle  \langle  \langle
	f^k_i,x^k_i
	\rangle :i<n_k \rangle ,h_k \right\rangle
	:k\in{\omega}
	\right\} $
	of pairs  $ \left\langle   \left\langle
	\left\langle f_i,x_i \right\rangle :i<n \right\rangle ,h \right\rangle $
	 satisfying
	\ref{df:large}.1--5.

	If $k\in{\omega}$ take
	\begin{displaymath}
	B^0_k=h^{-1}_k \left\{ 0 \right\} \cup \left\{ f^k_i({\nu}):i<n_k,{\nu}\in\dom(f
	^k_i)\mbox{ and }
	 \left\{ {\nu},x^k_i \right\} \notin E_{\beta} \right\}
	\end{displaymath}
	and
	\begin{displaymath}
	B^1_k=h^{-1}_k \left\{ 1 \right\} \cup \left\{ f^k_i({\nu}):i<n_k,{\nu}\in\dom(f
	^k_i)\mbox{ and }
	 \left\{ {\nu},x^k_i \right\} \in E_{\beta} \right\}.
	\end{displaymath}
	Applying  lemma \ref{lm:partition} ${\omega}$-many times we can find
	partitions $(C^0_k,C^k_1)$, $k<{\omega}$, of ${\omega}{\beta}$ such that
	taking
	\begin{displaymath}
	I^+_{\beta}
	     =\Cl\bigl((I_{\beta}\cup \left\{ C^1_k:k\in{\omega}\right\},
	            \bigcup_{{\gamma}\le{\beta}}\fcal_{\gamma}
	         \bigr)
	\end{displaymath}
	and
	\begin{displaymath}
	J^+_{\beta}
	     =\Cl\bigl((I_{\beta}\cup \left\{ C^0_k:k\in{\omega}\right\},
	            \bigcup_{{\gamma}\le{\beta}}\fcal_{\gamma}
	         \bigr)
	\end{displaymath}
	the set
	$ \left\{ S_{{\omega}{\gamma}}:{\gamma}\le{\beta}\right\}$
	is uncovered by  $I^+_{\beta}\cup J^+_{\beta}$.

	We can assume that $B^i_k{\subset} C^i_k$ for $i<2$ and $k<{\omega}$
	because $B^0_k\in J_{\beta}$ and $B^1_k\in I_{\beta}$.
	Take
	\begin{displaymath}
	E_{\alpha}=E_{\beta}\cup \left\{  \left\{ {\nu},{\omega}{\beta}+n \right\} :{\nu
	}<{\omega}{\beta},
	n\in{\omega}\mbox{ and }{\nu}\in B^1_n \right\}
	\end{displaymath}
	and
	\begin{displaymath}
	\fcal_{\alpha}=\emptyset.
	\end{displaymath}
	By the construction of
	$G_{\alpha}= \left\langle {\omega}{\alpha},E_{\alpha}\right\rangle$,
	it follows that
	${\omega}{\alpha}$ is large, so (II) holds. On the other hand
	\begin{displaymath}
	I_{\alpha}= \left\{ X\cup Y:X\in I^+_{\beta},Y\in[{\omega}{\alpha}]^{<{\omega}}
	\right\}
	\end{displaymath}
	and
	\begin{displaymath}
	J_{\alpha}= \left\{ X\cup Y:X\in J^+_{\beta},Y\in[{\omega}{\alpha}]^{<{\omega}}
	\right\} .
	\end{displaymath}
	so $ \left\{ S_{{\omega}{\gamma}}:{\gamma}<{\alpha} \right\} $ is uncovered by
	$I_{\alpha}\cup J_{\alpha}$.
	Finally $S_{{\omega}{\alpha}}$ is cofinal in ${\omega}{\alpha}$
	but the elements of $I_{\alpha}\cup J_{\alpha}$ are all bounded, so
	the induction hypothesis (I) also holds.
	\end{case}

	The construction is done. Take $E=\cup \left\{ E_{\alpha}:{\alpha}<\oo \right\}
	$
	and $G= \left\langle \oo,E \right\rangle $. By (I), $G$ is non-trivial.
	Finally, we must prove that $G$ is quasy smooth. Consider a set
	$Y{\subset}\oo$. The following lemma is almost trivial.
	\begin{lemma}
	\label{lm:split}
	For each ${\alpha}<\oo$ either $Y\cap({\omega}{\alpha}+{\omega})$
	or $({\omega}{\alpha}+{\omega}){\setminus} Y$ is large.
	\end{lemma}

	\proof
	Assume on the contrary that there are  pairs
	$ \left\langle  \left\langle  \left\langle f_i, x_i \right\rangle :i<n \right\rangle ,h \right\rangle $  and $ \left\langle  \left\langle  \left\langle f_i, x_i
	 \right\rangle :n\le i<n+k \right\rangle ,h' \right\rangle $
	showing that neither $Y\cap({\omega}{\alpha}+{\omega})$
	nor $({\omega}{\alpha}+{\omega}){\setminus} Y$ is large.
	Then
	$ \left\langle  \left\langle  \left\langle f_i, x_i \right\rangle :i<n+k \right\rangle ,h\cup h' \right\rangle $ shows that ${\omega}{\alpha}+{\omega}$
	is not large.
	\eproof

	So we can assume that the set
	\begin{displaymath}
	L= \left\{ {\alpha}<\oo:Y\cap({\omega}{\alpha}+{\omega})\mbox{
is large} \right\} 
	\end{displaymath}
	is uncountable and to complete the proof of theorem \ref{th:dplus}
	it is enough to show that in this case
	$G\cong G[Y]$.
	By $\dplus$, we can find a club subset $C{\subset} L'$
	such that $Y\cap {\omega}{\gamma}\in M_{\gamma}$,
	$C\cap{\omega}{\gamma}\in M_{\gamma}$ and ${\omega}{\gamma}={\gamma}$
	whenever ${\gamma}\in C$. We can assume that $0\in C$.

	Write $C= \left\{ {\gamma}_{\nu}:{\nu}<\oo \right\} $.
	By induction on ${\nu}<{{\omega}_1}$, we will construct functions
	$f_{\nu}$ such that
	\begin{enumerate}\alabel
	\item $f_{\nu}:G_{{\gamma}_{\nu}}\cong G_{{\gamma}_{\nu}}[Y]$,
	$f_{\nu}\in\fcal_{{\gamma}_{\nu}}$,
	\item $ \left\langle f_{\mu}:{\mu}<{\nu} \right\rangle
	\in M_{\sup \left\{ {\gamma}_{\mu}+1:{\mu}<{\nu}\right\}}$.
	\end{enumerate}
	Take $f_0=\emptyset$.
	If ${\nu}={\mu}+1$, then let
	$f_{\nu}=g^{Y\cap{\omega}{\gamma}_{\nu},f_{\mu}}$.
	If ${\nu}$ is limit, then put
	$f_{\nu}=\cup \left\{ f_{\mu}:{\mu}<{\nu} \right\} $.
	Clearly (a) and (b) remains valid.
	Finally put $f=\cup\left\{f_{\nu}:{\nu}<{{\omega}_1} \right\} $.
	Then $f:G\cong G[Y]$, so the theorem is proved.
	\eproof

	\section{A model without quasi-smooth graphs}
	\label{sc:no_sm}
	Given an Aronszajn-tree $T=\<\oo,\prec\>$ define the poset $\qcal_T$
	as follows:
	 the underlying set of $\qcal_T$ consists  of all  functions
	$f$ mapping a finite subset of $\oo$ to ${\omega}$ such that
	$f^{-1}\{n\}$ is antichain in $T$ for each $n\in{\omega}$.
	The ordering on $\qcal_T$ is as expected: $f\le_{\qcal_T}g$
	iff $f\supset g$.
	For ${\gamma}<\oo$ denote by $T_{\gamma}$ the set of elements of $T$
	with height ${\gamma}$.
	Take $T_{<{\delta}}=\bigcup\limits_{{\gamma}<{\delta}}T_{\gamma}$.
	If $x\in T_{{\delta}}$ and ${\gamma}<{\delta}$, let
	$x\rest {\gamma}$ be the unique element of $T_{\gamma}$ which is comparable
	with $x$.
	We write  $\ccal$ for the poset $\<\mbox{Fin}(\oo,2),\supset\>$, that is, forcing with
	$\ccal$   adds $\oo$-many Cohen reals to the ground model.
	\begin{theorem}
	\label{th:no_sm}
	If ZF is consistent, then so is ZFC + ``there are no non-trivial quasi-smooth
	graphs on $\oo$''.
	\end{theorem}

	\proof
	Assume that GCH holds in the ground model.  Consider a finite support
	iteration
	$\<P_i,Q_j:i\le\oot,j<\oot\>$
	satisfying (a)--(c) below:
	\begin{enumerate}\alabel
	\item If $j<\oot$ is even, then $Q_j=\ccal$.
	\item If $j<\oot$ is odd, then
	$V^{P_j}\models$ {\em``$Q_j=\qcal_{T_j}$ for some Aronszajn-tree $T_j$''}.
	\item
	$V^{P_\oot}\models${\em ``every Aronszajn tree is special''}.
	\end{enumerate}
	We will show that  $V^{P_\oot}$ does not
	contain non-trivial, quasi-smooth graphs on $\oo$.

	To start with we introduce some notation.
	Consider a graph  $G=\<V,E\>$. For $x\in V$ define the
	function
	$\typ(x):V\setm\{x\} \to 2$
	 by the equation
	$G(x)=\typ(x)^{-1}\{1\}$.
	Given $A\subs V$ write $\typ(x,A)=\typ(x)\rest A$.

	If $A\subs V$ and $t\in 2^A$,  take
	$\rel t;=\{x\in V\setm A:\typ(x,A)=t\}
	$
	and $\rls t;=\{x\in V\setm A:|\typ(x,A)\triangle t|<\omega\}$.
	For $x\in V$ and $A\subs V$ put
	$
	\twin x;A;=\rel \typ(x,A);
	$.

	For $A\subs V$ define the equivalence relation $\eqr A;$ on
	$V\setm A$ as follows:
	\begin{displaymath}
	x \eqr A; y \mbox{ \ iff \ }
	|\typ(x,A)\triangle \typ(y,A))|<\omega.
	\end{displaymath}

	For $x\in V\setm A$ denote by $\eqc x;A;$ the equivalence class of
	$x$ in $\eqr A;$.
	Clearly  $\eqc x;A;=\rls {\tp(x,A)};$.
	Write $\eqclass G;A;$ for the family of  equivalence
	classes of $\eqr A;$.

	We divide $\kcal$ into three subclasses, $\kcal_0$, $\kcal_1$ and
	$\kcal_2$,   and investigate them
	 separately to show that $V^{P_\oot}\models$
	``{\em $(\forall G\in \kcal_i)\ G$ is not quasi-smooth}'' for $i<3$.
	Take
	\begin{displaymath}
	\kcal_0=\{G\in\kcal:\exists A\in\br\oo;{\omega};
	\ |\eqclass G;A;|=\oo\},
	\end{displaymath}
	\begin{displaymath}
	\kcal_1=\{G\in\kcal:\forall A\in\br\oo;{\omega};
	\ \exists x\  |\oo\setm\eqc x;A;|<\oo\}
	\end{displaymath}
	and
	\begin{displaymath}
	\kcal_2=\kcal\setminus (\kcal_0\cup \kcal_1).
	\end{displaymath}

	\subsection{$G\in\kcal_0$}

	First we recall a definition of \cite{AS}.
	\begin{definition}
	A poset $P$ is {\em stable} if
	\begin{displaymath}
	\forall B\in\br P;{\omega};
	\ \exists B^*\in \br P;{\omega};
	\ \forall p\in P\ \exists p'\le p
	\ \exists p^*\in B^*
	\ \forall b\in B\  (p'\compp b\mbox{ iff }p^*\compp b).
	\end{displaymath}
	We will
	say that $p'$ and $p^*$ are {\em twins for $B$} and that {\em $B^*$
	shows the stability of $P$ for $B$}.
	\end{definition}

	\begin{lemma}
	\label{lm:pstbl}
	$P_{\oot}$ is stable.
	\end{lemma}

	\proof
	First let us remark that it is enough to prove that
	both  $\ccal$ and $\qcal_T$ are stable for any Aronszajn-tree for
	in \cite{AS} it was proved that  any finite support
	iteration of stable, c.c.c. posets is  stable.

	It is clear that $\ccal$ is stable. Assume that
	$T$ is an Aronszajn tree and   $B\subs\br\qcal_T;{\omega};$.
	Fix a countable ordinal ${\delta}$
	with $\{\dom(p):p\in B\}\subs T_{<{\delta}}$ and take
	$B^*=\{p\in\qcal_T:\dom(p)\subs T_{< {\delta}+{\omega}}\}$.
	It is not hard to see that $B^*$ shows the stability of $P$ for $B$.
	\eproof

	For  $G\in\kcal$ take
	$G\in\kcal^*_0$ iff
	there is an $A\in\br\oo;{\omega};$
	such that
	 the set
	$\{x:|\eqc x;A;|\le\omega\}$ is uncountable.

	Given $G\in\kcal_0$ we will write
	$G\in\kcal'_0$ iff
	there are disjoint sets $A_0, A_1\in\br\oo;{\omega};$
	such that
	\begin{enumerate}\arablabel
	\item
	$x \eqr A_0;y$ iff  $x\eqr A_1;y$ for each  $x,y\in\oo\setm  A_0\cup A_1$,
	\item the set
	$\{x:|\eqc x;A_0;|\le\omega\}$ is uncountable.
	\end{enumerate}

	\begin{lemma}
	\label{lm:kp0}
	Assume CH. If $G\in\kcal'_0$, then there is a partition
	$(V_0,V_1)$ of $\oo$ so that for each stable c.c.c.
	poset $P$ we have
	\begin{displaymath}
	V^P\models \mbox{``$G$ is not isomorphic to  $G[V_i]$ for $i\in 2$''}.
	\end{displaymath}
	\end{lemma}

	\proof
	Pick $A_0,A_1\in\br\oo;{\omega};$ witnessing $G\in\kcal'_0$.
	Write $A=A_0\cup A_1$.
	Take $E=E(G)$.

	Let $\k$ be a large enough regular cardinal and fix
	an increasing sequence $\ooseq N;\n;$ of countable, elementary
	submodels of $\hcal_{\k}$ such that
	\begin{enumerate}\rlabel
	\item $G,A,A_0,A_1\in N_0$,
	\item $\seq N;\n;\m;\in N_{\mu}$ for ${\mu}<\oo$.
	\end{enumerate}

	For $x\in\oo\setm A$ take
	\begin{displaymath}
	\rank(x)=\min\{\n:x\in N_{\nu}\}.
	\end{displaymath}

	 Fix a partition $(S_0,S_1)$ of $\oo$ with $|S_0|=|S_1|=\oo$.
	Take
	\begin{displaymath}
	V_i=A_i\cup\rank^{-1}S_i
	\end{displaymath}
	for $i\in 2$.

	We show that the partition $(V_0,V_1)$ works.

	Assume on the contrary that $P$ is a stable c.c.c. poset, $\dot f$ is a
	$P$-name of a function, $p_0\in P$ and
	\begin{displaymath}
	p_0\force \mbox{``}\dot f:G\cong G[V_0]\mbox{''}.
	\end{displaymath}

	Without loss of generality we can assume that $p_0=1_P$.
	Now for each $c\in A_0$
	choose  a maximal antichain  $J_c\subs P$ and a
	 function  $h_c:J_c\to V$ such that
	 $q\force$``$\dot f^{-1}(\hat c)=\widehat{h_c(r)}$'' for each $q\in J_c$.

	Take
	$
	B=\bigcup\{J_c: c\in A_0\}
	$
	and pick a countable $B^*\subs P$ showing the stability of $P$ for $B$.

	For $b\in P$ define the partial function
	$\dtb:\oo\to 2^{A_0}$ as follows.
	Let $x\in\oo$. If there is a function
	$t\in 2^{A_0}$ so that
	\begin{enumerate}\alabel
	\item $t(c)=1$
	$\Longleftrightarrow$ for each $q\in I_c$ if $q$ and $b$ are compatible
	conditions, then $\{x,h_c(q)\}\in E$,
	\item $t(c)=0$
	$\Longleftrightarrow$ for each $q\in I_c$ if $q$ and $b$ are compatible
	conditions, then $\{x,h_c(q)\}\notin E$,
	\end{enumerate}
	then take $\dtb (x)=t$. Otherwise $x\notin\dom\, \dtb$.

	\begin{sublemma}
	\label{lm:dtb}
	If $p\force$``$\fdot(x)=y$'', then there are $p'\le p$ and $b\in B^*$
	such that $b$ and $p'$ are twins for $B$ and
	$\dtb(x)=\tp(y,A_0)$.
	\end{sublemma}

	\proof
	By the choice of $B^*$, we can find a $p'\le p$ and a $b\in B^*$
	so that $p'$ and $b$ are twins for $B$.
	Let $c\in A_0$.
	For each $q\in J_c$, if $q$ and $p'$ are compatible in $P$, then
	$\{y,c\}\in E$  iff $\{x,h_c(q)\}\in E$ , because,
	taking  $r$ as a common extension of $q$ and $p'$, we have
	{\em $r\force $``$\dot f(\hat x)=\hat y\mbox{ and }
	\dot f(\widehat{h_q(c)})=\hat c$''.}
	So
	$\{y,c\}\in E$ iff for each $q\in I_c$ if $q$ and $p'$ are compatible,
	then $\{x,h_c(q)\}\in E$.
	But $p'$ and $b$ are twins for $\bigcup\{J_c:c\in A_0\}$,
	so $\dtb(x)=\dt p;(x)=\tp(y,A_0)$.
	\eproof

	\begin{sublemma}\label{lm:dtlarge}
	There is a $b\in B^*$ such that
	$$\tagg{*}
	|\{t\in\ran\, \dtb:|\rls t;|\le\omega \}|=\oo.
	$$
	\end{sublemma}

	\proof
	Let $\gcal$ be a $P$-generic filter over $V$.
	Put
	$$\fcal=\{\tp(y,A_0):y\in V_0\setm A_0, |\eqc y;A_0;|\le\omega\}.$$
	Then $|\fcal|=\oo$, so we can write
	$\fcal=\ooset t;\nu;$. Fix sequences
	$\ooseq p;\nu;\subs\gcal$, $\ooseq x;\nu;\subs\oo$ and
	$\ooseq y;\nu;\subs\oo$ such that
	$p_{\nu}\force$``$\fdot(x_{\nu})=y_{\nu}$'' and
	$\tp(y_{\nu},A_0)=t_{\nu}$.
	By sublemma \ref{lm:dtb},
	$$
	\bigcup\limits_{b\in B^*}\ran\, \dtb\supseteq \fcal.
	$$
	But $B^*$ is countable, so we can find a $b\in B^*$ satisfying
	$(*)$ above.
	\eproof

	Fix $b\in B^*$ with property $(*)$.
	Consider the structure
	\begin{displaymath}
	\ncal=\<P\rest (B\cup B^*),B,B^*,\<J_c,h_c:c\in A_0\>\>.
	\end{displaymath}
	By  CH, there is a ${\nu}<\oo$ with $\ncal\in N_{\nu}$.
	Pick ${\mu}\in S_1\setm \n$.
	Since $G$, $\ncal$, $b\in N_{\nu}$, it follows that
	$\dtb\in N_{\nu}\subs N_{\mu}$.
	By $(*)$ and (ii), there is a
	$$
	t\in\ran\, \dtb\cap (N_{\mu}\setm
	\bigcup\limits_{\xi<\mu}N_{\xi})
	$$
	with $|\rls t;|\le\omega$.
	Then
	$$\tagg{\dag}
	\rls t;\subs N_{\mu}\setm\bigcup\limits_{\xi<\mu}N_{\xi}.
	$$
	Pick $x\in\oo$ with $\dtb(x)=t$.
	Find $p\in\gcal$ and $y\in V_0$ such that
	 $p\le b$ and $p\force$``$\fdot(x)=y$''.
	By sublemma \ref{lm:dtb}, there are $p'\le p$ and
	$b'\in B^*$ such that $p'$ and $b'$ are twins for $B$
	and $\dt b';(x)=\tp(y,A_0)$. But $p\le b$, so
	$\dtb(x)=\dt b';(x)$. Indeed, let $c\in A_0$ and
	assume that $\dt b';(x)(c)=1$. Pick $q\in J_q$ which is compatible with
	$b'$. By the definition of $\dt b';$, it follows that
	$\{h_c(q),x\}\in E$. Since $p'$ and $b'$ are twins for $B$, so
	$p'$ and $q$ also have a common extension $q'$ in $P$. But
	$p'\le p\le b$, so $q'$ witnesses that $b$ and $q$ are compatible.
	Thus, by the definition of $\dtb$, we have $\dtb(x)(c)=1$.

	Thus $\tp(y,A_0)=\dt b';(x)=\dtb(x)=t$.
	By $(\dag)$, this implies that $\rank(y)=\mu$.
	But, by the construction of the partition $(V_0,V_1)$,
	there are no $y\in V_0$ with $\rank(y)=\mu$.
	Contradiction, the lemma is proved.
	\eproof

	\begin{lemma}
	\label{lm:kkprime}
	Assume CH. If $G\in\kcal^*_0$, then
	\newline
	$V^{\ccal}\models${\em ``
	there is a partition
	$(V_0,V_1)$ of $\oo$ so that
	\newline
	\makebox[4cm]{}for each stable c.c.c.
	poset $P$ we have:}
	\begin{displaymath}
	V^{\ccal*P}\models \mbox{``$G$ is not isomorphic to  $G[V_i]$ for $i\in 2$''.\ ''}
	\end{displaymath}
	\end{lemma}

	\proof . Fix a set
	$A\in\br\oo;{\omega};$ witnessing $G\in\kcal^*_0$
	 and a bijection $f:A\to\omega$  in V.
	Let $r:\omega\to 2$ be the characteristic function of a Cohen real
	 from   $V^{\ccal}$.
	Take $A_i=(f\circ r)^{-1}\{i\}$ for $i<2$. Then  $(A_0,A_1)$
	is a partition of $A$. Using a trivial density argument
	we can see that $x\noteqr A;y$ implies $x\noteqr A_i;y$ for $i<2$
	and for $x,y\in\oo\setm A$.
	Thus
	 $V^{\ccal}\models$ ``$A_0$ and $A_1$
	witness $G\in\kcal'_0$''.
	Applying lemma \ref{lm:kp0} in $V^{\ccal}$ we
	get the desired partition of $\oo$.
	\eproof

	\begin{lemma}\label{lm:k0ks}
	In $V^{P_{\omega_2}}$, if $G\in\kcal_0$ is quasi-smooth, then
	$G\in\kcal^*_0$.
	\end{lemma}

	\proof
	Choose a set
	$A\in\br\oo;{\omega};$ witnessing $G\in\kcal_0$
	 and a bijection $f:A\to\omega$.
	Pick $\alpha<\oot$, $\alpha$ is even, with
	$A$, $f$, $G\in V^{P_{\alpha}}$. From now on we work in
	$V^{P_{\alpha}}$. Let $\{\eqc x_{\nu};A;:\nu<\oo\}$ be an enumeration of
	the equivalence classes of $\eqr A;$.
	Fix a partition $(I_0,I_1)$  of $\oo$ into uncountable pieces.
	Let $r:\omega\to 2$ be the characteristic function of a Cohen real
	 from   $V^{P_{\alpha}*\ccal}$.
	Take $A_i=(f\circ r)^{-1}\{i\}$ for $i<2$. Then  $(A_0,A_1)$
	is a partition of $A$. Using a trivial density argument
	we can see that $x\noteqr A;y$ implies $x\noteqr A_i;y$ for $i<2$
	and for $x,y\in\oo\setm A$.
	For $i\in 2$ put
	$$
	B_i=A_i\cup\{x_{\nu}:\nu\in I_i\}\cup
	\{\eqc x_{\nu};A;\setm\{x_{\nu}\}:\nu\in I_{1-i}\}.
	$$
	Clearly $(B_0,B_1)$ is a partition of $\oo$ and
	$$B_i\cap \eqc x_\nu;A_i;=B_i\cap \eqc x_\nu;A;=\{x_\nu\}.$$
	So $G[B_i]\in\kcal^*$. But $G$ is quasi-smooth, so
	$G\cong G[B_i]$ for some $i\in 2$ in
	$V^{P_{\oot}}$. Thus $G\in\kcal^*_0$ is proved.
	\eproof


	\subsection{$G\in\kcal_1$}

	We say that a poset $P$ has property \Prr\
	iff for each sequence $\ooseq p;\n;\subs P$
	there exist disjoint sets $U_0,U_1\in\br\oo;\oo;$ such that whenever
	$ {\alpha}\in U_0$ and
	${\beta}\in U_1$ we have $p_{\alpha}\|_{{}_P}p_{\beta}$.

	\begin{lemma} $\ccal$ has property \Prr.
	\end{lemma}

	Indeed, $\ccal$ has property K.

	\begin{lemma} If $T$ is an Aronszajn-tree, then $\qcal_T$
	has property \Prr.
	\end{lemma}

	\proof
	Let $\ooseq p;{\alpha};\subs P$ be given.
	We can assume that there are  a stationary set $S\subs \oo$,
	 $p^*\in \qcal_T$,  ${\gamma}^*<\oo$, $n\in{\omega}$
	and $\{z_i:i<n\}\subs T$  such that for each ${\alpha}\in S$
	\begin{enumerate}\alabel
	\item $x\rest{\alpha}\in\dom(p_{\alpha})$ for each
	$x\in \dom(p_{\alpha})$ with $\height_{T}(x)\ge\alpha$,
	\item $p_{\alpha}\rest T_{<{\alpha}}=p^*$,
	\item $|\dom(p_{\alpha})\cap T_{\alpha}|=n$,
	\item writing
	$\dom(p_{\alpha})\cap T_{\alpha}=\{x^{\alpha}_0,\ldots,x^{\alpha}_{n-1}\}$,
	$x^{\alpha}_0<_{\mbox{\tiny On}}\ldots<_{\mbox{\tiny On}}x^{\alpha}_{n-1}$,
	the sequence $\<p_{\alpha}(x^{\alpha}_0),\ldots,p_{\alpha}(x^{\alpha}_{n-1})\>$
	is independent from ${\alpha}$,
	\item ${\gamma}^*<{\alpha}$ and
	the elements $x^{\alpha}_0\rest {\gamma}^*,\ldots,x^{\alpha}_{n-1}\rest {\gamma}
	^*$ are pairwise distinct,
	\item $x^{\alpha}_i\rest {\gamma}^*=z_i$ for $i<n$.
	\end{enumerate}

	For each  ${\beta}<\oo$ and
	${\bar y}=\<y_0,\ldots,y_{n-1}\>\in (T_{\beta})^n$
	take
	\begin{displaymath}
	S_{\bar y}=
	\{{\alpha}\in S\setm {\beta}: x^{\alpha}_i\rest {\beta}=y_i
	\mbox{ for each $i<n$}\}.
	\end{displaymath}

	Let
	\begin{displaymath}
	C^*=\{\dd<\oo\setm{\gamma}^*:\forall {\beta}<\dd\ \forall
	\bar y\in (T_{\beta})^n\ (|S_{\bar y}|\le{\omega}\to S_{\bar y}\subs\dd)\}.
	\end{displaymath}
	Now $\{p_{\alpha}:{\alpha}\in S\cap C^*\}$ are $\oo$ members of $P$, so
	for some ${\alpha}<{\beta}\in S\cap C^*$ the conditions  $p_{\alpha}$ and
	$p_{\beta}$ are compatible. Since
	$p_{\alpha}(x^{\alpha}_l)=p_{\beta}(x^{\beta}_l)$, $x^{\alpha}_l$ and $x^{\beta}
	_l$
	are incomparable in $T$ for $l<n$.
	So for some $\n<{\alpha}$,
	$x^{\alpha}_l\rest{\nu}\ne x^{\beta}_l\rest{\nu}$
	whenever  $l<n$.
	On the other hand, for $l\ne m <n$ we have
	$x^{\alpha}_l\rest{\nu}\ne x^{\beta}_m\rest{\nu}$
	because
	$x^{\alpha}_l\rest{\gamma}^*=z_l\ne z_m=x^{\beta}_m\rest{\gamma}^*$.
	Take $y^a_l=x^{\alpha}_l\rest {\nu}$ and
	 $y^b_l=x^{\beta}_l\rest {\nu}$ for $l<n$ and write
	$\bar a=\<y^a_0,\ldots,y^a_{n-1}\>$,
	$\bar b=\<y^b_0,\ldots,y^b_{n-1}\>$ .
	The elements $\{y^a_i,y^b_i:i<n\}$ are pairwise different,
	so for each ${\alpha}'\in S_{\abar}$ and ${\beta}'\in S_{\bbar}$
	the conditions $p_{{\alpha}'}$ and $p_{{\beta}'}$ are compatible.
	But
	 $|S_{\bar a}|=|S_{\bar b}|=\oo$, because
	${\alpha}\in S_{\bar a}$, ${\beta}\in S_{\bar b}$, ${\nu}<{\alpha}$ and
	${\alpha}\in C^*$.
	 \eproof

	A poset $P$ is called {\em well-met} if any two compatible elements $p_0$
	and $p_1$
	of $P$ have a greatest lower bound denoted by $p_0\land p_1$.

	\begin{lemma}
	\label{lm:prod}
	Assume that the poset $P$ has property \Prr\
	 and $V^P\models$ ``{\em the poset $Q$ has property \Prr}''.
	Let $\{\<p_{\alpha},q_{\alpha}\>:{\alpha}<\oo\}\subs P*Q$.
	Then there are disjoint sets $U_0,U_1\in\br\oo;\oo;$ such that
	for each ${\gamma}\in U_0$ and ${\delta}\in U_1$ the conditions
	$\<p_{\gamma},q_{\gamma}\>$ and $\<p_{\delta},q_{\delta}\>$ are compatible,
	in other words, $p_{\gamma}$ and $p_{\delta}$ have a common extension
	 $p_{{\gamma},{\delta}}$ in $P$  with
	$p_{{\gamma},{\delta}}\force$``$q_{\gamma}\parallel_Q q_{\delta}$''.
	If $P$ is well-met, then we can find
	 conditions $\{p'_{\alpha}:{\alpha}\in U_0\cup U_1\}$  in $P$
	with
	$p'_{\alpha}\le p_{\alpha}$ such that
	$p'_{\gamma}\land p'_{\delta}\force$``$q_{\gamma}\parallel_Q q_{\delta}$''
	for each ${\gamma}\in U_0$ and ${\delta}\in U_1$.
	\end{lemma}

	\proof
	Let $\dot U$ be a $P$-name for the set
	$U=\{{\alpha}:p_{\alpha}\in \gcal_P\}$, where $\gcal_P$ is the $P$-generic
	filter.
	Since $P$ satisfies c.c.c., there is a $p^*\in P$ with
	$p^*\force\mbox{``$|\dot U|=\oo$''}$.
	Since $V^P\models${\em ``$Q$ has property \Prr''},
	there is a condition $p\le p^*$ and there are $P$-names such that
	$p\force$
	{\em ``$V_i=\{\adot^i_{\gamma}:{\gamma}<\oo\}\in\br U;\oo;$,
	for $i\in 2$,
	and $q_{\adot^0_{\gamma}}$ and $q_{\adot^1_{\delta}}$
	are compatible whenever  ${\gamma},{\delta}\in\oo$''}.
	Choose conditions $p^*_{\gamma}\le p$ and ordinals
	${\beta}^0_{\gamma}$, ${\beta}^1_{\gamma}$,
	with  $p^*_{\gamma}\force$``
	$\adot^i_{\gamma}={\hat\beta}^i_{\gamma}$'' for $i<2$.

	Now consider the sequence $A=\{p^*_{\gamma}:{\gamma}<\oo\}$. Since
	$P$ has property \Prr, there are disjoint, uncountable sets
	$C_0,C_1\subs A$ such that $p^*_{\gamma}$ and $p^*_{\delta}$
	are compatible  whenever ${\gamma}\in C_0$ and
	${\delta}\in C_1$. Take
	$U_i=\{{\beta}^i_{\gamma}:{\gamma}\in C_i\}$ for  $i\in 2$.
	We can assume that $U_0\cap U_1=\empt$.
	Let ${\gamma}\in C_0$ and ${\delta}\in C_1$
	 and let  $p''$ be a common extension of $p^*_{\gamma}$ and
	$p^*_{\delta}$. Then $p''\force$``{\em
	${\beta}^0_{\gamma},{\beta}^1_{\delta}\in \dot U$, that is,
	$p_{{\beta}^0_{\gamma}}$
	 and $p_{{\beta}^1_{\delta}}$ are in $\gcal_P$}'', so $p''$,
	must be a common extension of $p_{{\beta}^0_{\gamma}}$
	 and $p_{{\beta}^1_{\delta}}$.
	So  $p''\force$``${\beta}^0_{\gamma}\in V_0$ and
	${\beta}^1_{\delta}\in V_1$'', thus
	$p''\force$``$q_{{\beta}^0_{\gamma}}$ and
$q_{{\beta}^1_{\delta}}$ are compatible in $Q$'', 
	so $\<p_{{\beta}^0_{\gamma}},q_{{\beta}^0_{\gamma}}\>
	\parallel_{P*Q}\<p_{{\beta}^1_{\delta}},q_{{\beta}^1_{\delta}}\>$.

	Suppose   that  $P$ is well-met.
	Take $p'_{{\beta}^0_{\gamma}}=p^*_{\gamma}\land p_{{\beta}^0_{\gamma}}$
	and
	$p'_{{\beta}^1_{\delta}}=p^*_{\delta}\land p_{{\beta}^1_{\delta}}$.
	It works because we can use $p^*_{\gamma}\land p^*_{\delta}$ as $p''$ in the argument of
	the previous paragraph.
	\eproof

	\begin{lemma}
	\label{lm:prit}
	If $\<R_{\alpha}:{\alpha}\le{\mu},S_{\beta}:{\beta}<{\mu}\>$
	is a finite support iteration such that
	$V^{R_{\alpha}}\models$``$S_{\alpha}$ has property \Prr''
	for ${\alpha}<{\mu}$, then $R_{\mu}$ has property \Prr, as well.
	\end{lemma}

	\proof
	We prove this lemma by induction on $\mu$. The successor  case  is
	covered by lemma \ref{lm:prod}. Assume that $\mu$ is limit.
	Let $\ooseq p;\xi;\subs R_{\mu}$. Without loss of generality
	we can assume that $\<\supp(p_{\xi}):\xi<\oo\>$ forms a $\Delta$-system
	with kernel $d$. Fix $\nu<\mu$ with $d\subs \nu$. By the induction
	hypothesis, the poset $R_\nu$ has property \Prr, so there exist
	disjoint sets $U_0$, $U_1\in\br\oo;\oo;$ such that
	whenever
	$ {\xi}\in U_0$ and
	${\eta}\in U_1$ we have $p_{\xi}\|_{{}_{R_\nu}}p_{\eta}$.
	But $p_{\xi}\|_{{}_{R_\nu}}p_{\eta}$ implies
	$p_{\xi}\|_{{}_{R_\mu}}p_{\eta}$ because
	$\supp(p_{\xi})\cap\supp(p_{\eta})\subs\nu$, so $R_\mu$ has property \Prr,
	as well.
	\eproof

	The previous lemmas yield the following corollary.
	\begin{lemma}
	\label{lm:proo}
	$P_\oot$ has property \Prr.
	\end{lemma}

	Given $G=\<\oo,E\>\in\kcal_1$ and $\xi$, $\alpha$, $\beta\in\oo$ with
	$\xi\in\alpha\cap\beta$ take
	\[
	\dabb G;\xi;\alpha;\beta;=
	\{{\nu}\in\xi:
	\{{\alpha},{\nu}\}\in E\mbox{ iff }\{{\beta},{\nu}\}\notin E\}.
	\]

	\begin{lemma}\label{lm:dabb}
	If $G\in\kcal_1$, then
	$$\tagg{*}
	\forall\xi\in\oo\ \exists\eps\xi;\in\oo\
	\forall\alpha,\beta\in\oo\setm\eps\xi;\
	|\dabb G;\xi;\alpha;\beta;|<\omega.
	$$
	\end{lemma}

	\proof Since $G\in\kcal_1$, we have an $x\in\oo$ with
	$|\oo\setm [x]_\xi|<\oo$. Choose $\eps\xi;\in\oo\setm\xi$
	with $\oo\setm [x]_\xi\subs\eps\xi;$. It works because
	$\alpha$, $\beta>\eps \xi;$ implies $\alpha$, $\beta\in [x]_{\xi}$.
	\eproof

	The bipartite graph $\<\oo\times 2,
	\{\{\<{\nu},0\>,\<{\mu},1\>\}:{\nu}<{\mu}<\oo\}\>$
	 will be denoted by $\halfg$.

	\begin{lemma}
	\label{lm:half}
	If $G\in\kcal_1$, then neither $G$ nor its complement may have
	a --- not necessarily spanned --- subgraph isomorphic to
	$\halfg$.
	\end{lemma}

	\proof
	Let $G=\<\oo,E\>$. Write $E(\alpha)=\{\xi\in\oo:\{\xi,\alpha\}\in E\}$.
	Assume on the contrary that $A$, $B\in\br\oo;\oo;$ are disjoint sets
	such that $\{\alpha,\beta\}\in E$ whenever
	$\alpha\in A$ and $\beta\in B$ with
	$\alpha<\beta$. Without loss of generality we can assume that
	$(A\setm \alpha+1)\cap \epsilon(\alpha)=\empt$ for each $\alpha\in A$.
	Write $A=\ooset \alpha;\xi;$. Then for $\xi\in\oo$ the set
	$F(\xi)=(A\cap \alpha_{\xi})\setm E(\alpha_{\xi+1})$
	is  finite because $\alpha_{\xi+1}>\epsilon(\alpha_{\xi})$ and
	$(A\cap \alpha_{\xi})\setm E(\beta)=\empt$ for all but countable many
	$\beta\in B$. By Fodor's lemma, we can assume that $F(\xi)=F$
	for each $\xi\in S$, where $S$ is a stationary subset of $\oo$
	containing limit ordinals only.
	Let $T=\{\xi\in S:F\subs\alpha_\xi\}$ and take
	$W=\{\alpha_{\xi+1}:\xi\in T\}$.
	Then  $G[W]$ is an uncountable complete subgraph of $G$.
	Contradiction.
	\eproof

	\begin{lemma}
	\label{lm:k1}
	If $G\in\kcal_1$ and $V^{\ccal}\models$``{\em $Q$  has property \Prr}'',
	then
	\begin{displaymath}
	V^{\ccal* Q}\models \mbox{``}G\not\cong G[f^{-1}\{i\}]\mbox{ for $i<2$'',}
	\end{displaymath}
	where $f:\oo\to 2$ is the $\ccal$-generic function over $V$.
	\end{lemma}

	\proof
	Assume on the contrary that
	\begin{displaymath}
	\<p,q\>\force\mbox{``$\hhbar:G\cong G[f^{-1}\{0\}]$''}.
	\end{displaymath}
	To simplify our notations, we will write
	$E$ for $E(G)$,  $\dab \xi;\alpha;\beta;$ for $\dabb G; \xi;\alpha;\beta;$
	and $\ep \xi;$ for $\eps\xi;$.


	Let $C_0=\{{\delta}<\oo:\xi<\delta\mbox{ implies } \ep\xi;<\delta\}$.
	Clearly $C_0$ is club. Take $C_1=\{{\delta}<\oo:\<p,q\>\force
	\mbox{``}\hhbar''{\hat\delta}=f^{-1}\{0\}\cap{\hat\delta}\mbox{''}\}$.
	Since $\ccal*Q$ satisfies c.c.c, the set $C_1$
	 is club.
	Put $C_2=C_0\cap C_1$.

	Now for each ${\alpha}<\oo$ let
	$\delta_\alpha=\min (C_2\setm \alpha+1)$ and choose a condition
	$\<p_{\alpha},q_{\alpha}\>\le\<p,q\>$ and a countable ordinal
	${\gamma}_{\alpha}$
	such that
	\begin{displaymath}
	\<p_{\alpha},q_{\alpha}\>\force\mbox{``}\hhbar(\hat\delta_\alpha)=
	{\hat\gamma}_{\alpha}\mbox{''}.
	\end{displaymath}

	Since $\gamma_\alpha\ge\delta_\alpha>\ep \alpha;$ for each
	$\alpha\in\oo$, we can fix a stationary set
	$S\subs \oo$ and a finite set $D$ such that
	 $\dab{\alpha};{\delta_{\alpha}};{\gamma}_{\alpha};=\mbox{D}$ for each
	${\alpha}\in S$.
	Since $\ccal$ is well-met, applying lemma \ref{lm:prod}
	we can find disjoint uncountable subsets $S_0,S_1\subs S$
	and a sequence $\<p'_{\alpha}:{\alpha}\in S_0\cup S_1\>\subs\ccal$  with
	$p'_{\alpha}\le p_{\alpha}$ such that
	$p'_{\alpha}\land p'_{\beta}\force$ ``$q_{\alpha}\parallel_Q q_{\beta}$''
	for  each ${\alpha}\in S_0$
	and ${\beta}\in S_1$.

	We can assume that
	the sets $\{\dom(p'_{\alpha}):{\alpha}\in S_0\}$ and
	$\{\dom(p'_{\beta}):{\beta}\in S_1\}$ form  $\Delta$-systems with
	kernels $d_0$ and $d_1$,
	respectively.

	Take
	$
	Y^0_{\xi}=\{{\alpha}\in S_0: \{{\xi},\delta_{\alpha}\}\in E\}
	$
	and
	$Y^1_{\xi}=\{{\alpha}\in S_1: \{{\xi},\delta_{\alpha}\}\notin E\}$
	for  ${\xi}<\oo$.
	Write
	$
	Y_i=\{{\xi}<\oo:|Y^i_{\xi}|=\oo\}
	$
	and $Z_i=\oo\setm Y_i$ for $i<2$.

	By \ref{lm:half}, the sets $Z_i$ are countable.
	Pick ${\xi}\in C_2$ with $D\cup d_0\cup d_1\cup Z_0\cup Z_1\subs{\xi}$.
	Let ${\xi}'=\min (C_2\setm {\xi}+1)$ and
	${\xi}''=\min (C_2\setm {\xi'}+1)$ .
	Since  $d_0\cup d_1\subs\xi$
	and $|Y^0_{\xi}|=|Y^1_{\xi}|=\oo$,
	we can choose ${\alpha}_i\in Y^i_{\xi}\setm {\xi}''$
	with $\dom(p'_{{\alpha}_i})\cap [{\xi},{\xi}')=\empt$ for $i=0,1$.
	The set
	$W=\dab{{\xi^\prime}};\delta_{{\alpha}_0};\delta_{{\alpha}_1};
	\cap [{\xi},{\xi}')$is finite because
	$\delta_{{\alpha}_i}\ge\alpha_i\ge\xi''>\ep\xi';$ for $i<2$.
	Choose a $\ccal$-name $q$ such that
	$p'_{{\alpha}_0}\land p'_{{\alpha}_1}\force$
	{\em``$q$ is a common extension of $q_{{\alpha}_0}$ and $q_{{\alpha}_1}$
	in $Q$''} and
	take
	\begin{displaymath}
	r=\<p'_{{\alpha}_0}\cup p'_{{\alpha}_1}\cup\{\<{\nu},1\>:{\nu}\in W\},q\>.
	\end{displaymath}

	Since $W\cap(\dom(p'_{{\alpha}_0})\cup\dom(p'_{{\alpha}_1}))=\empt$,
	$r$ is a condition.

	Pick  a condition $r'\le r$ from $\ccal*Q$
	and an  ordinal ${\eta}$ such that
	$r'\force $``$\hhbar({\hat\xi})={\hat\eta}$''.
	Now  $\eta\in [\xi,\xi')$ because $\xi$, $\xi'\in C_1$.
	Since
	\begin{displaymath}
	r'\force
	\mbox{``}\hhbar(\delta_{{\alpha}_i})=\hat\gamma_{{\alpha}_i},
	\hhbar({\hat\xi})={\hat\eta}\mbox{ and }
	\hhbar\mbox{ is an isomorphism''},
	\end{displaymath}
	so $\{\delta_{{\alpha}_0},\xi\}\in E$ and $\{\delta_{{\alpha}_1},\xi\}\notin E$
	imply that
	$\{{\gamma}_{{\alpha}_0},\eta\}\in E$ and
	$\{{\gamma}_{{\alpha}_1},\eta\}\notin E$.
	But $D_{\alpha_i}(\delta_{{\alpha}_i},{\gamma}_{{\alpha}_i})=D$ and
	$D\subs\xi$
	so $\{\delta_{{\alpha}_0},\eta\}\in E$ and
	$\{\delta_{{\alpha}_1},\eta\}\notin E$,
	that is, $\eta\in W$.
	But $r\force $``{\em $\ran(\hhbar)=f^{-1}\{0\}$ and
	$f^{-1}\{0\}\cap\hat W=\empt$}'', contradiction.
	\eproof

	\subsection{$G\in\kcal_2$}

	Given a non-trivial graph  $G=\<V,E\>$ with $V\in\br\oo;\oo;$ define
	\begin{displaymath}
	\Gamma(G)=\{\dd\in\oo:\exists {\alpha}\in V\ {\alpha}\ge\dd\mbox{ and }
	|\mbox{twin}_G({\alpha},V\cap\dd)|\le{\omega}\}.
	\end{displaymath}
	 The following lemma obviously holds.
	\begin{lemma}
	\label{lm:inv}
	If $G_0$ and $G_1$ are graphs on uncountable subsets of $\oo$,
	$G_0\cong G_1$, then
	$\Gamma(G_0)=\Gamma(G_1)\mbox{\rm\ mod }\nstat_{\oo}$.
	\end{lemma}

	\begin{lemma}
	\label{lm:part} Given $G\in\kcal\setm\kcal_0$ and
	 $S\subs\oo$ there is a partition $(V_0,V_1)$ of $\oo$ such that
	$\Gamma(G[V_0])\subs S\mod \nstat_{\oo}$ and
	$\Gamma(G[V_1])\subs \oo\setm S\mod \nstat_{\oo}$.
	\end{lemma}

	\proof
	Let $\k$ be a large enough regular cardinal and fix
	an increasing, continuous sequence $\ooseq N;\n;$ of countable, elementary
	submodels of $\hcal_{\k}=\<H_{\k},\in\>$ such that
	$G,S\in N_0$ and $\<N_{\nu}:{\nu}\le{\mu}\>\in N_{{\mu}+1}$ for ${\mu}<\oo$.
	Write ${\gamma}_{\nu}=N_{\nu}\cap\oo$ and $C=\ooset{\gamma};{\nu};$.
	Take $V_0=\bigcup\limits_{{\nu}\in S}({\gamma}_{{\nu}+1}\setm {\gamma}_{\nu})$
	and $V_1=\oo\setm V_0=
	\bigcup\limits_{{\nu}\in \oo\setm S}({\gamma}_{{\nu}+1}\setm {\gamma}_{\nu})$.

	It is enough to prove that $\Gamma(G[V_0])\subs S\mod\nstat_{\oo}$.
	Assume that
	${\gamma}_{\nu}\in\Gamma(G[V_0])$,
	${\gamma}_{\nu}={\nu}$, ${\alpha}\ge{\gamma}_{\nu}$, ${\alpha}\in V_0$
	and $|\mbox{twin}_{G[V_0]}({\alpha},{\gamma}_{\nu}\cap V_0)|={\omega}$.
	Since $G$, ${\nu}$, ${\gamma}_{\nu}\cap V_0\in N_{{\nu}+1}$ and
	$|\eqclass G;V_0\cap {\gamma}_{\nu};|\le{\omega}$,
	we have $\typ_{G[V_0]}({\alpha},{\gamma}_{\nu}\cap V_0)\in N_{{\nu}+1}$
	and so
	$\mbox{twin}_{G[V_0]}({\alpha},{\gamma}_{\nu})\subs N_{{\nu}+1}$ as well.
	Thus ${\alpha}\in {\gamma}_{{\nu}+1}\setm {\gamma}_{\nu}$.
	Hence ${\alpha}\in V_0$ implies ${\gamma}_n={\nu}\in S$
	which was to be proved.
	\eproof

	\begin{lemma}
	\label{lm:k2}
	If $G\in\kcal\setm \kcal_0$ and $\Gamma(G)\ne\empt \mod\nstat_{\oo}$,
	then $G$ is not quasi-smooth.
	\end{lemma}

	\proof
	Assume that $S=\Gamma(G)$ is stationary and let
	$(S_0,S_1)$ be a partition of $S$ into stationary subsets.
	By lemma \ref{lm:part}, there is a partition $(V_0,V_1)$ of $\oo$
	 with $\Gamma(G([V_i])\cap S\subs S_i$.
	Then $G[V_i]$ and $G$ can not be isomorphic by lemma \ref{lm:inv}.
	\eproof

	Let us remark that  $G\in\kcal_2$ iff $G\in\kcal\setm \kcal_0$ and
	there is an $A\in\br\oo;{\omega};$ and $x\in\oo\setm A$
	such that
	 $|\eqc x;A;|=|\oo\setm \eqc x;A;|=\oo$.

	Given $G\in\kcal_2$ we will write $G\in\kcal'_2$ iff
	there are two disjoint, countable subsets of $\oo$, $A_0$ and $A_1$,
	 and there is an $x\in\oo$,
	such that
	$|\eqc x;{A_0};|=|\oo\setm \eqc x;{A_0};|=\oo$ and
	$\eqc x;{A_0};\setm A_1= \eqc x;{A_1};\setm A_0$.

	\begin{lemma}
	\label{lm:k2k2p}
	If $G\in\kcal_2$, then $G\in(\kcal'_2)^{V^{\ccal}}$.
	\end{lemma}

	\proof
	Assume that $A\in\br\oo;{\omega};$ and $x\in\oo$ witness
	$G\in\kcal_2$ in the ground model.
	Fix a bijection $f:A\to\omega$  in V.
	Let $r:\omega\to 2$ be the characteristic function of a Cohen real
	from   $V^{\ccal}$.
	Take $A_i=(f\circ r)^{-1}\{i\}$.
	By a simple density argument,
	we can see that
	$\eqc x;A_0;=\eqc x;A;=\eqc x;A_1;$.
	Thus
	$A_0$,  $A_1$ and $x$
	show that $g\in\kcal'_2$.
	\eproof

	\begin{lemma}
	\label{lm:gml}
	Assume that every Aronszajn tree is special.
	If $G\in\kcal'_2$, then
	there is a partition $(V_0,V_1)$ of  $\oo$ such that
	$\Gamma(G[V_i])$  is stationary for $i<2$.
	\end{lemma}

	\proof
	Choose $A_0,A_1$ and $x$ witnessing $G\in\kcal'_2$.
	Let $A=A_0\cup A_1$.
	Take $C_0=\eqc x;A_0;\setm A$,  $C_1=(\oo\setm \eqc x;A_0;)\setm A$
	and consider the partition trees $\tcal_i$ of $G[C_i]$
	for $i\in 2$ (see definition
	\ref{df:ptree}).
	These trees are
	 Aronszajn-trees because $G$ is non-trivial.
	Fix  functions $h_i\colon C_i\to{\omega}$ specializing $\tcal_i$.
	We can find natural numbers $n_0$ and $n_1$ such that
	 the sets
$S_i=\{\nu:h_i^{-1}\{n_i\}\cap(\tcal_i)_{\nu}\ne\empt\}$ are
stationary, that is, $h_i^{-1}\{n_i\}$ meets stationary many level of  
	$\tcal_i$.
	Take $B_i=h^{-1}_i\{n_i\}$
	and $Y_i=\{c\in C_i:\exists b\in B_i\ c\preceq_{\tcal_i}b\}$.

	Pick any ${\delta}\in S_i$. Let $b\in B_i\cap (\tcal_i)_{\delta}$.
	If $c\in Y_i\setm (\tcal_i)_{<{\delta}}$, $c\ne b$,
	then $c\rest {\delta}\ne b$ by the construction of $Y_i$.
	So
	$\typ_{G[Y_i]}(c,(\tcal_i)_{<{\delta}})=
	\typ_{G[Y_i]}(c\rest {\delta},(\tcal_i)_{<{\delta}})\ne
	\typ_{G[Y_i]}(b,(\tcal_i)_{<{\delta}})$
	by the definition of the partition tree.
	This means that
	$\mbox{twin}_{G[Y_i]}(b,(\tcal_i)_{<{\delta}})=\{b\}$. Thus
	${\delta}\in\Gamma(G[Y_i])$ provided $(\tcal_i)_{<{\delta}}\subs {\delta}$
	and $b\ge{\delta}$. But these
	requirements exclude only a non-stationary
	subset of $S_i$. So $\Gamma(G[Y_i])\supset S_i\mod \nstat_{\oo}$.

	Let  $V_i= Y_i\cup A_i\cup (C_{1-i}\setm Y_{1-i})$ for $i\in 2$ and
	consider the partition $(V_0,V_1)$ of $\oo$ .
	If $z\in V_i\setm (Y_i\cup A_i)$, then $\typ(z,A_i)\ne \typ(b,A_i)$
	for any $b\in B_i$ because
	$C_0\subs \eqc x; A_i;$ and $C_1\subs \oo\setm \eqc x;A_i;$. So
	$\Gamma(G[V_i])\supset S_i\mod\nstat_{\oo}$ holds.
	\eproof


	Now we are ready to conclude the proof of theorem \ref{th:no_sm}.
	We will work in  $V^{P_{\oot}}$.
	Assume that $G\in\kcal$. We must show that $G$ is not quasi-smooth.

	Pick a ${\nu}<\oot$
	with $G\in(\kcal)^{V^{P_{\nu}}}$ and $Q_{\nu}=\ccal$.
	Assume first that $G\in(\kcal_0)^{V^{P_{\nu}}}$.
	If $G$ were quasi-smooth in $V^{P_{\oot}}$,
	$G\in (\kcal^*_0)^{P_{\oot}}$ would hold by lemma \ref{lm:k0ks}.
	So we can assume that
	$G\in (\kcal^*_0)^{P_{\nu}}$.
	Since $P_{\oot}$ is a stable, c.c.c. poset, so is
	$P_{\oot}/P_{{\nu+1}}$.
	So, by lemma \ref{lm:kkprime}, there is a partition $(V_O,V_1)$
	of $\oo$ in $V^{P_{\nu+1}}$ such that
	$V^{P_{\oot}}
	\models$``{\em $G$ is not isomorphic to $G[V_i]$ for $i<2$}''.

	Assume  that  $G\in(\kcal_1)^{V^{P_{\nu}}}$.
	Since $P_{\oot}$ has property \Prr, so is
	$P_{\oot}/P_{{\nu+1}}$.
	Thus, by lemma \ref{lm:k1}, the partition $(V_O,V_1)$
	of $\oo$ given by the $Q_{\nu}$-generic Cohen reals  in $V^{P_{\nu+1}}$
	has the property that
	$V^{P_{\oot}}
	\models$``{\em $G$ is not isomorphic to $G[V_i]$ for $i<2$}''.

	Finally assume that $G\in(\kcal_2)^{V^{P_{\nu}}}$. By lemma
	\ref{lm:k2k2p}, we have $G\in(\kcal'_2)^{V^{P_{\nu+1}}}$. Since
	$P_{\oot}$ satisfies c.c.c, it follows that $G\in(\kcal'_2)^{V^{P_{\oot}}}$.
	So applying lemma \ref{lm:gml} we can find a partition $(V_0,V_1)$ of
	$\oo$ such that both $\Gamma(G[V_0])$ and $\Gamma(G[V_1])$  are
	stationary. Thus, by lemma \ref{lm:k2},  neither $G[V_0]$ nor
	$G[V_1]$ are quasi-smooth.   So $G$ itself can not be quasi-smooth.
	\eproof

	\end{document}